\pdfoutput=1
\documentclass[a4paper]{amsart}
\usepackage{amsaddr} 
\usepackage{xspace} 
\usepackage{ifthen}
\usepackage[hidelinks]{hyperref}
\usepackage{cleveref}
\usepackage{graphicx}
\usepackage{float}
\usepackage{geometry}

\newgeometry{vmargin={25mm,18mm}, hmargin={1in,1in}}     
\newcommand*{\M}[1]{\ensuremath{#1}\xspace} 

\newcommand*{\x}{\times}
\newcommand*{\mi}[1]{\mathbf{#1}} 
\newcommand*{\st}[1]{\mathbb{#1}} 
\newcommand*{\rv}[1]{\mathsf{#1}} 
\newcommand*{\te}[2][]{\left\lbrack{#2}\right\rbrack_{#1}}
\newcommand*{\tte}[2][]{\lbrack{#2}\rbrack_{#1}}

\newcommand*{\diag}[2][]{\left\langle{#2}\right\rangle_{#1}}
\newcommand*{\prob}[2]{\M{\mathbb{P}\!\left\lbrack\left.{#1}\;\right\vert\;{#2}\right\rbrack}} 
\newcommand*{\deq}{\M{\mathrel{\mathop:}=}} 
\newcommand*{\deqr}{\M{=\mathrel{\mathop:}}} 
\newcommand{\T}[1]{\text{#1}} 
\newcommand*{\QT}[2][]{\M{\quad\T{#2}\ifthenelse{\equal{#1}{}}{\quad}{#1}}} 
\newcommand*{\ev}[3][]{\mathbb{E}_{#3}^{#1}\!\left\lbrack{#2}\right\rbrack}
\newcommand*{\evt}[3][]{\mathbb{E}_{#3}^{#1}#2}
\newcommand*{\cov}[3][]{\ifthenelse{\equal{#1}{}}{\mathbb{V}_{#3}\!\left\lbrack{#2}\right\rbrack}{\mathbb{V}_{#3}\!\left\lbrack{#2,#1}\right\rbrack}}
\newcommand*{\dev}[3][]{\ifthenelse{\equal{#1}{}}{\mathbb{D}_{#3}\!\left\lbrack{#2}\right\rbrack}{\mathbb{D}_{#3}\!\left\lbrack{#2,#1}\right\rbrack}}
\newcommand*{\covt}[2]{\mathbb{V}_{#2}{#1}}
\newcommand*{\devt}[2]{\mathbb{D}_{#2}{#1}}
\newcommand*{\gauss}[2]{\mathcal{N}\!\left({#1,#2}\right)}
\newcommand*{\uni}[2]{\mathcal{U}\!\left({#1,#2}\right)}

\newcommand*{\modulus}[1]{\M{\left\lvert{#1}\right\rvert}}

\newcommand*{\set}[1]{\M{\left\lbrace{#1}\right\rbrace}} 
\newcommand*{\setbuilder}[2]{\M{\left\lbrace{#1}\: \big\vert \:{#2}\right\rbrace}}
\newcommand*{\uniti}{\lbrack 0,1\rbrack}
\newcommand*{\unitia}{\lbrack -0.001,1.001\rbrack}

\DeclareMathOperator*{\ish}{ish}
\DeclareMathOperator*{\sob}{sob}
\DeclareMathOperator*{\oak}{oak}

\begin{document}

\title{Sobol' Matrices For Multi-Output Models With Quantified Uncertainty}
\author[R.A. Milton and S.F. Brown]{Robert A. Milton and Solomon F. Brown}
\email{r.a.milton@sheffield.ac.uk, s.f.brown@sheffield.ac.uk}
\address{Department of Chemical, Materials and Biological Engineering, The University of Sheffield, Sheffield~S1~3JD, United Kingdom}
\thanks{This work was funded by the EPSRC under grant reference EP/V051458/1.}
\begin{abstract}
    Variance based global sensitivity analysis measures the relevance of inputs to a single output using Sobol' indices. This paper extends the definition in a natural way to multiple outputs, directly measuring the relevance of inputs to the linkages between outputs in a correlation-like matrix of indices. The usual Sobol' indices constitute the diagonal of this matrix. Existence, uniqueness and uncertainty quantification are established by developing the indices from a putative multi-output model with quantified uncertainty. Sobol' matrices and their standard errors are related to the moments of the multi-output model, to enable calculation. These are benchmarked numerically against test functions (with added noise) whose Sobol' matrices are calculated analytically.
  \end{abstract}
\maketitle

\section{Introduction}\label{sec:Intro}
    This paper is concerned with analyzing the results of experiments or computer simulations in a design matrix of $M\geq 1$ input axes (columns) and $L\geq 1$ output axes (columns) over $N$ samples (rows). Global Sensitivity Analysis (GSA) \cite{Razavi2021} examines the relevance of the various inputs to the various outputs. When pursued via ANOVA decomposition of a single output, this leads naturally to the well known Sobol' indices, which have by now been applied across most fields of science and engineering \cite{Saltelli2019,Ghanem2017}. This paper extends the definition in a natural way to multiple outputs $L\geq 1$. 

    Sobol' indices apportion the variance of a scalar output to the influence of various inputs. There is no generally accepted extension of Sobol' indices to multi-outputs $L>1$. 
    The Sobol' matrices developed in the present work apportion the correlation matrix of outputs -- rather than the variance of a single output -- to the influence of various inputs. The diagonal elements are identically the classical Sobol' indices, expressing the relevance of inputs to the output variables themselves. The off-diagonal elements express the relevance of inputs to the correlation between outputs. This may be of vital interest when outputs are, for example, yield and purity of a product, cost and efficacy of a process, or perhaps a single output measured at various times. The Sobol matrices reveal (amongst other things) which inputs it is worthwhile varying in an effort to alter the linkage between outputs.
    
    Prior work on Sobol' indices with multiple outputs \cite{Gamboa.etal2013,GarciaCabrejo2014,Xiao2017,Cheng2019} has settled ultimately on just the diagonal elements of the covariance matrix, so covariance between outputs remains unexamined. Output covariance has been tackled in previous studies by performing principal component analysis (PCA) on outputs prior to performing GSA on the resulting diagonal covariance matrix of outputs \cite{Campbell2006}. This has been used in particular to study synthetic multi-outputs which are actually the dynamic response of a single output over time \cite{Lamboni2011,Zhang2020}. Hierarchical sensitivity analysis \cite{Xu2020} employs prior mapping to independent subsystems. In all these cases, output covariance is essentially removed prior to GSA. Only very recently has GSA attempted to incorporate output covariance at all \cite{Liu2021}, as an ingredient in a total fluctuation. The present work is the first to attempt GSA directly on output covariance itself. This direct approach is deliberately and explicitly motivated by the straightforward desire to examine the influence of various inputs on the correlation between outputs.

    Accurate calculation of Sobol' indices for a single output is computationally expensive and may require 10,000+ samples \cite{Lamoureux.etal2014, Oakley.OHagan2004}. A (potentially) more efficient approach is calculation via a surrogate model, such as a Gaussian Process (GP) \cite{Sacks.etal1989, Rasmussen.Williams2005}, Polynomial Chaos Expansion (PCE) \cite{Ghanem.Spanos1997,Xiu.Karniadakis2002,Xiu2010}, low-rank tensor approximation \cite{Chevreuil.etal2015,Konakli.Sudret2016}, or support vector regression \cite{Cortes.Vapnik1995,Cheng2019}. Semi-analytic expressions for Sobol' indices are available for scalar PCEs \cite{Sudret2008}, and the diagonal elements of multi-output PCEs \cite{GarciaCabrejo2014}.
    Semi-analytic expressions for Sobol' indices of GPs have been provided in integral form by \cite{Oakley.OHagan2004} and alternatively by \cite{Chen.etal2005}. These approaches are implemented, examined and compared in \cite{Marrel.etal2009,Srivastava.etal2017}. Both \cite{Oakley.OHagan2004,Marrel.etal2009} estimate the errors on Sobol' indices in semi-analytic, integral form. Fully analytic, closed form expressions have been derived without error estimates for uniformly distributed inputs \cite{Wu.etal2016a} using a GP with a radial basis function kernel.

    In this paper we develop Sobol' matrices and their error estimates in a general setting, then relate these to the moments of a surrogate model. This is to facilitate efficient calculation with limited data, when some form of surrogate is really indispensable. The resulting analytic formulae have been implemented for a multi-output Gaussian process (MOGP), allowing us to benchmark implementation against test functions with known Sobol' matrices. This is important, because practitioners are rarely privileged with an exact and reliable exogenous model linking outputs to inputs. Rather, GSA is a preliminary stage in reducing empirical observations to a tractable, analyzable representation. For example, it is frequently used to guide model reduction and the design of experiments, as well as or as part of optimization and risk assessment. Under these circumstances, the GSA itself must be robust and reliable, providing well benchmarked measures of its own uncertainty and domain of validity.

    The contents of this paper are as follows. The next Section introduces multi-output models with quantified uncertainty (MQUs) as a context which embraces any explicit or black-box function (simulation, experiment), surrogate (emulator, response surface, meta-model) or regression which supplies a probability distribution on its outputs. This encompasses multi-output models with no uncertainty (MNUs) as MQUs with zero variance.
    \Cref{sec:GSAI,sec:GSAM} define Sobol' indices for MQUs, then extend them to Sobol' matrices. We seek to emphasise that this work is the natural -- arguably unique -- extension of the \emph{definition} of a Sobol' index. Otherwise it may seem that the classical Sobol' indices occupying the diagonal could be extended to a full correlation matrix in a multitude of arbitrary ways. \Cref{sec:GSAS,sec:Mom} derive mean and standard error estimates on the Sobol' matrices of an MQU, then expresses them via the moments of the MQU. The latter is the key step towards implementation. The perspective throughout is essentially Bayesian, regarding uncertainty as undetermined or unobserved input(s).

    \Cref{sec:Func} introduces an $L=9$ output MNU composed of known test functions, and examines its Sobol' matrices. Test MQUs consisting of the test MNU plus Gaussian noise are used to benchmark a GP implementation of Sobol' matrices in \Cref{sec:Benchmarks}. This highlights some advantages and limitations of the GP implementation. Conclusions are drawn in \Cref{sec:Conc}.

\section{Multi-output models with quantified uncertainty}\label{sec:MQU}
    The purpose of this Section is to set the scene for this study, serve as a glossary of notation, and prove the formal foundation of the construction to follow. Our notation is not quite standard, but is preferred for lightness and fluency once grasped. Regarding notation and other topics, the intention is to facilitate computation by combining efficient, tensorized operations such as {\tt einsum} \cite{Numpy2025,PyTorch2025,Tensorflow2025} into numerically stable calculations. As such, the devices employed herein should be familiar to practitioners of machine learning.

    The preference throughout this paper is for uppercase constants and lowercase variables.
    Square bracketed quantities such as $\tte[\mi{M\x N}]{u}$ are tensors subscripted by (the cartesian product of) their ranks, for bookkeeping. Each rank is expressed as an ordinal of axes, a tuple which is conveniently also a naive set and may be subtracted from others as such
    \begin{equation*} \label[definition]{def:MQU:m}
        \begin{aligned}
            \mi{0} \deq () \subseteq{} \mi{m}&\deq (0,\ldots ,m-1) \subseteq \mi{M} \deq (0,\ldots ,M-1) \\   
            \mi{M-m} &\deq \setbuilder{m^{\prime} \in \mi{M}}{m^{\prime} \notin \mi{m}} = (m, \ldots, M-1)
        \end{aligned}
    \end{equation*}
    GSA decomposes input space by families of axes, and this notation facilitates that. It does assume that input axes are expediently ordered already, which might impair convenience but not generality. Conventional indexing of components is permitted under this scheme, by allowing singleton axes to be written in italic without parentheses.

    From a Bayesian perspective a determined tensor is a function of determined input(s), a random variable (RV) a function of undetermined input, and a stochastic process (SP) a function of both determined and undetermined inputs. In any case, all sources of (sans-serif) uncertainty may be gathered in a single input $\tte[M]{\rv{u}}$, final to the definition of input space:
    \begin{equation} \label[definition]{def:MQU:input}
        \begin{aligned}
            \mi{M} &\deq (0,\ldots,M-1) \QT{determined inputs on the unit interval} &\te[\mi{M}]{u}\in \uniti^{M} \\
            M &\deqr \mi{M+1-M} \quad\QT{an undetermined input on the unit interval} &\te[M]{\rv{u}} \sim \uni{0}{1}
        \end{aligned}   
    \end{equation}
    Exponentiation is always categorical -- repeated cartesian $\x$ or tensor $\otimes$ -- unless stated otherwise.
    It is crucial to this work that all $\mi{M+1}$ inputs vary completely independently of each other -- they are in no way codependent or correlated.
    Uncertainty is distributed uniformly $\tte[M]{\rv{u}} \sim \uni{0}{1}$ to exploit the ``universality of the uniform'' \cite[pp.224]{Blitzstein2019}, better known as the inverse probability integral transform \cite{Angus1994}. This explicitly provides a bijective mapping to any RV $\tte[M]{\rv{x}}$ by inverting its cumulative distribution function (CDF) $\mathrm{P} \colon \tte[M]{\rv{x}} \mapsto \tte[M]{\rv{u}}$ in the generalized sense
    \begin{equation*}
        \mathrm{P}^{-1}(\tte[M]{u}) \deq \inf \setbuilder{\tte[M]{x}}{\mathrm{P}(\tte[M]{x}) \geq \tte[M]{u}}        
    \end{equation*}
    This mapping $\mathrm{P}^{-1} \colon \tte[M]{\rv{u}} \mapsto \tte[M]{\rv{x}}$ is the quantile function of $\tte[M]{\rv{x}}$ on $\uniti$, and embedding it in an output model $y$ can express any continuous RV at a given determined input $\tte[\mi{M}]{u}$. In other words, the output model could be any continuous SP.
    An amenity of \Cref{def:MQU:input} is that each input $\te[m]{\rv{u}}$ is automatically a probability measure, for any $m\in\mi{M}$.

    Throughout this work, expectations are subscripted by the input axes marginalized
    \begin{equation*}
        \ev{\bullet }{\#} \deq \int_{\te[\#]{0}}^{\te[\#]{1}} \te{\bullet} \, \mathrm{d}\te[\#]{u} \qquad \forall \# \subseteq \mi{M+1} \T{ or } \# \in \mi{M+1}
    \end{equation*}
    Covariances invoke the tensor product (summing the ranks of the arguments), and carry the subscript of their underlying expectation
    \begin{equation*}
        \cov[\te{*}]{\te{\bullet}}{\#} \deq \ev{\te{\bullet}\otimes\te{*}}{\#} - \ev{\bullet}{\#} \otimes \ev{*}{\#} \qquad \forall \# \subseteq \mi{M+1} \T{ or } \# \in \mi{M+1}
    \end{equation*}
    The covariance of anything with itself is expressed with a single argument $\cov{\bullet}{\mi{\#}} \deq \cov[\bullet]{\bullet}{\mi{\#}}$, as is customary.

    A multi-output model with quantified uncertainty (MQU) is defined as any Lebesgue integrable function on input space \Cref{def:MQU:input} obeying
    \begin{equation} \label[definition]{def:MQU:y}
        \begin{gathered}
            y \colon \uniti^{M+1} \rightarrow \st{R}^{L} \QT{such that} \\ \covt{\evt{\te[l]{y(\te[\mi{M+1}]{\rv{u}})}}{M}}{\mi{M}} > \evt{\evt{\te[l]{y(\te[\mi{M+1}]{\rv{u}})}}{M}}{\mi{M}} = 0 \quad \forall l\in\mi{L}            
        \end{gathered}
    \end{equation}
    This forces every output component $l\in\mi{L}$ to depend on at least one determined input: constant functions and pure noise are not allowed, to avoid division by zero in GSA.
    Without loss of generality we have also offset the output  $\tte[\mi{L}]{y}$ to have mean $\tte[\mi{L}]{0}$. This is formally unnecessary, but it can be crucial to the numerical accuracy of computations.

    Our notation will reflect machine learning practice \cite{Numpy2025,PyTorch2025,Tensorflow2025} which facilitates parallel computation by tensorizing function application over batch dimensions such as $\mi{N}$ according to
    \begin{equation*}
        \te[\mi{L}\x n]{y(\te[\mi{M+1\x N}]{\rv{u}})} \deq y(\te[\mi{M+1}\x n]{\rv{u}}) \qquad \forall n \in \mi{N}
    \end{equation*}

    To elucidate the meaning and purpose of an MQU, and set the stage for GSA, consider a design matrix of $\tte[\mi{M\x N}]{u}$ inputs alongside $\tte[\mi{L\x N}]{Y}$ outputs:
    \begin{equation*}
        \begin{array}{lllll}
            \te[0\x 0]{u} & \cdots & \te[M-1\x 0]{u} \\
            \phantom{t}\vdots & \ddots & \phantom{t}\vdots \\
            \te[0\x N-1]{u} & \cdots & \te[M-1\x N-1]{u} \\
        \end{array}
        \begin{array}{lllll}
            \te[0\x 0]{Y} & \cdots & \te[L-1\x 0]{Y} \\
            \phantom{t}\vdots & \ddots & \phantom{t}\vdots \\
            \te[0\x N-1]{Y} & \cdots & \te[L-1\x N-1]{Y} \\
        \end{array}
    \end{equation*}
    The sole purpose of an MQU is actually to provide a single number for any such design matrix of $N\in\st{Z}^{+}$ samples (rows): namely, the probability that the output samples result from the input samples, to within any given uncertainty $\te[\mi{L\x N}]{e} \geq \te[\mi{L\x N}]{0}$
    \begin{equation*}
        \prob{\te[\mi{L\x N}]{E} \geq {y(\te[\mi{M+1\x N}]{\rv{u}})-\te[\mi{L\x N}]{Y}} \geq -\te[\mi{L\x N}]{E}}{\te[\mi{M\x N}]{\rv{u}}=\te[\mi{M\x N}]{u}}
    \end{equation*}
    This is what it means to generate $L$ quantifiably uncertain outputs from $M$ determined inputs.
    Following the Kolmogorov extension theorem \cite[pp.124]{Rogers.Williams2000}, the meaning of an MQU is nothing other than an SP. This includes zero variance SPs -- fully determined multi-output models of no uncertainty (MNUs) -- as MQUs $y(\tte[\mi{M+1}]{\rv{u}})$ which do not depend on the final, undetermined input $\tte[M]{\rv{u}}$. Kolmogorov's extension theorem formally identifies an SP -- a collection $y(\tte[\mi{M+1}]{\rv{u}}) = \setbuilder{y(\tte[\mi{M+1}]{\rv{u}})}{\tte[\mi{M}]{\rv{u}}=\tte[\mi{M}]{u}}$ of RVs indexed by determined inputs $\tte[\mi{M}]{u}$ -- with all its finite dimensional distributions $\setbuilder{y(\tte[\mi{M+1\x N}]{\rv{u}})}{\tte[\mi{M\x N}]{\rv{u}}=\tte[\mi{M\x N}]{u}}$. The latter viewpoint is the formal version of the random field \cite{Khoshnevisan2002} or random function interpretation of SPs \cite[pp.42]{Skorokhod2005} frequently alluded to in machine learning \cite{Rasmussen.Williams2005}.
    It can be a helpful perspective on the objects appearing in this work.

    In conception MQUs include any simulation, surrogate or regression which supplies a probability distribution (even a zero variance one) on its outputs. In practice, such an MQU is usually inferred from a design of experiments (DoE) which samples the determined inputs $\tte[\mi{M}]{u}$ from $\uni{0}{1}^{M}$. This is not restrictive, a preferred set of determined inputs $\tte[\mi{M}]{x}$ may have any sampling distribution whatsoever, so long as it is continuous and its $\mi{M}$ axes mutually independent. Simply take the CDF $\tte[\mi{M}]{\rv{u}} = \mathrm{P}(\tte[\mi{M}]{\rv{x}})$ before beginning, and apply the quantile function $\tte[\mi{M}]{\rv{x}}=\mathrm{P}^{-1}(\tte[\mi{M}]{\rv{u}})$ described above and in \cite{Angus1994} after finishing. In general and in summary, simply treat any $y^{\dag}(\tte[\mi{M+1}]{\rv{x}})$ as the pullback of an MQU by the CDF of $\tte[\mi{M+1}]{\rv{x}}$.

    To close this Section we shall formally secure the construction which follows.
    An MQU $y$ is Lebesgue integrable by \Cref{def:MQU:y}, so it must be measurable. All measures throughout this work are finite (probability measures, in fact), so integrability of $y$ implies integrability of $y^n$ for all $n \in \st{Z}^{+}$ \cite{Villani1985}.
    Therefore, Fubini's theorem \cite[pp.77]{Williams1991} allows all expectations to be taken in any order, which will be crucial later. We can thus safely construct any SP which is polynomial in $y$ and its marginals; and freely extract finite dimensional distributions from it by the Kolmogorov extension theorem \cite[pp.124]{Rogers.Williams2000}. This guarantees existence and uniqueness of every device in this paper.

\section{Sobol' indices}\label{sec:GSAI}
    This Section recapitulates the definition of classical Sobol' indices in the context of an MQU, to prepare the way for easy generalization to Sobol' matrices in the next Section. Sobol' indices are classically constructed \cite{Sobol2001} from a High Dimensional Model Representation (HDMR or Hoeffding-Sobol' decomposition \cite{Chastaing2011}). This is enlightening, but not strictly necessary. Instead, we shall swiftly construct the indices from an MQU. Later, this will allow us to derive standard errors on the Sobol' indices.

    In essence, GSA is performed by marginalizing (obscuring) determined input axes selectively. 
    A reduced model (mQU) is an SP defined as
    \begin{equation}\label[definition]{def:GSAI:y_m}
    \te[\mi{L}]{\rv{y}_{\mi{m}}} \deq \ev{y(\te[\mi{M+1}]{\rv{u}})}{\mi{M-m}}_{\mi{L}} \qquad \forall \mi{m} \subseteq \mi{M}
    \end{equation}
    When $\mi{M-m}=\mi{0}$, nothing is marginalized and $\rv{y}_{\mi{M}}$ is called the full model (MQU).
    The undetermined input $\tte[M]{\rv{u}}$ is never marginalized, and each mQU depends on $\mi{m}$ determined inputs $\tte[\mi{m}]{u}$ alongside $\tte[M]{\rv{u}}$. Any distributional assumption about inputs simply pulls back (factors through) expectations here with no impact whatsoever: $\uni{0}{1}$ is no more than a proxy for the CDF of any continuous distribution.

    An mQU $\te[\mi{L}]{\rv{y}_{\mi{m}}}$ is \emph{not} automatically converted to a new MQU over $\te[\mi{m+1}]{\rv{u}}$ by re-indexing $\te[M]{\rv{u}}$ to $\te[m]{\rv{u}}$. This would \emph{not} be the same as simply hiding columns $\mi{M-m}$ in the design matrix of the previous Section. In general the new MQU would \emph{not} adequately reflect the original full model, because any output variation due to marginalized inputs $\te[\mi{M-m}]{u}$ has simply been lost: it cannot be re-allocated to the new undetermined input $\te[m]{\rv{u}}$. This gets to the heart of GSA: assessing the degree to which a reduced model mQU mimics its full model MQU. Which is a matter of assessing how much output variation has been lost through marginalization.
    
    The marginal variance of an mQU is a tensor-valued RV (function of $\te[M]{\rv{u}}$) defined as
    \begin{equation}\label[definition]{def:GSAI:V_m}
        \cov{\rv{y}_{\mi{m}}}{\mi{m}}_{\mi{L\x L}} \deq \cov{\te[\mi{L}]{\rv{y}_{\mi{m}}}}{\mi{m}} 
    \end{equation}
    It is -- for all realizations $\te[M]{\rv{u}}$=$\te[M]{u}$ -- a symmetric tensor by definition and positive semi-definite by Jensen's inequality. We may therefore define a vector RV which is the square root of its diagonal, namely the standard deviation $\dev{\rv{y}_{\mi{m}}}{\mi{m}}_{\mi{L}}$
    \begin{equation}\label[definition]{def:GSAI:D_m}
        \dev{\rv{y}_{\mi{m}}}{\mi{m}}_{l} \deq \sqrt{\cov{\rv{y}_{\mi{m}}}{\mi{m}}_{lxl}} \quad > 0 \qquad \forall l \in \mi{L}
    \end{equation}
    This is positive definite due to the positive variance clause in the \Cref{def:MQU:y} of an MQU.

    The closed Sobol' index of scalar output $\tte[l]{y(\tte[\mi{M+1}]{\rv{u}})}$ with respect to input axes $\mi{m}$ is defined as an RV
    \begin{equation}\label[definition]{def:GSAI:S_m}
            \te[l]{\rv{S}_\mi{m}} \deq \frac{\cov{\rv{y}_{\mi{m}}}{\mi{m}}_{l\x l}}{\cov{\rv{y}_{\mi{M}}}{\mi{M}}_{l\x l}}
    \end{equation}
    The definition originated \cite{Sobol1993} for fully determined functions ($y$ independent of $\tte[M]{\rv{u}})$. It was later extended to random (i.e. undetermined) inputs \cite{Sobol.Kucherenko2005} and GPs \cite{Oakley.OHagan2004,Marrel.etal2009}. The latter work introduced the Sobol' index as an RV constructed from an important category of SP.

    The complement of a closed index is called a total index \cite{Homma1996}
    \begin{equation}\label[definition]{def:GSAI:ST_m}
        \te[l]{\rv{S}^{T}_\mi{M-m}} \deq 1 - \te[l]{\rv{S}_\mi{m}}
    \end{equation}

    Sobol' indices are readily interpreted, given that codependent inputs (reviewed in \cite{PWiederkehrThesis}) are beyond the scope of this work. By \Cref{def:GSAI:S_m} a closed index $S_{\mi{m}} \in \uniti$ measures the proportion of output variance captured by the reduced model $\rv{y}_{\mi{m}}$. This, in turn, indicates the influence or relevance of input axes $\mi{m}$. The full model $\rv{y}_\mi{M}$ explains everything explicable, so its Sobol' index is $\tte[l]{\rv{S}_{\mi{M}}}=\tte[l]{\rv{S}^{T}_{\mi{M}}}=1$. The void model $\rv{y}_\mi{0}$ is just the mean output -- an RV depending only on undetermined noise $\te[M]{\rv{u}}$ -- which explains nothing, so its Sobol' index is $\tte[l]{\rv{S}_{\mi{0}}}=\tte[l]{\rv{S}^{T}_{\mi{0}}}=0$.

    The closed index with respect to a single input axis $m\in \mi{M}$ is called a first-order index
    \begin{equation}\label[definition]{def:GSAI:S_mp}
        \te[l]{\rv{S}_{m}} \deq \te[l]{\rv{S}_\mi{m+1-m}}
    \end{equation}
    Because inputs may cooperate to affect the output, a closed index often exceeds the sum of its first-order contributions, obeying (for any realization $\tte[M]{\rv{u}}=\tte[M]{u}$)
    \begin{equation}\label{eq:GSAI:order}
        \sum_{\dot{m} \in \mi{m}} \te[l]{\rv{S}_{\dot{m}}} \leq \te[l]{\rv{S}_{\mi{m}}} \leq  \te[l]{\rv{S}^{T}_{\mi{m}}}
    \end{equation}
    The final inequality observes that a closed index only includes cooperation of input axes $\mi{m}$ with each other, whereas a total index also includes cooperation between $\mi{m}$ and $\mi{M-m}$ (but excludes cooperation of input axes $\mi{M-m}$ with each other). This is the only difference between closed and total indices, but it is an important one.

    Regarding the covariance between any reduced model and its full model, Fubini's theorem and the law of iterated expectations safeguard the intuitive identity
    \begin{equation} \label[definition]{def:GSAI:V_mM}
        \begin{aligned}
            \cov[\rv{y}_{\mi{M}}]{\rv{y}_{\mi{m}}}{\mi{M}}_{\mi{L\x L}} &\deq \ev{\te[\mi{L}]{\rv{y}_{\mi{m}}} \otimes \te[\mi{L}]{\rv{y}_{\mi{M}}}}{\mi{M}} - \evt{\te[\mi{L}]{\rv{y}_{\mi{m}}}}{\mi{M}} \otimes \evt{\te[\mi{L}]{\rv{y}_{\mi{M}}}}{\mi{M}} \\
            &\phantom{:}= \ev{\te[\mi{L}]{\rv{y}_{\mi{m}}} \otimes \te[\mi{L}]{\rv{y}_{\mi{m}}}}{\mi{m}} - \te[\mi{L}]{\rv{y}_{\mi{0}}} \otimes \te[\mi{L}]{\rv{y}_{\mi{0}}} \\
            &\deqr \cov{\rv{y}_{\mi{m}}}{\mi{m}}_{\mi{L\x L}}
        \end{aligned}
    \end{equation}
    Using this, the correlation between (the scalar output predicted by) an mQU and its parent MQU is an RV
    \begin{equation} \label[definition]{def:GSAI:R_m}
        \te[l]{\rv{R}_{\mi{mM}}} \deq \frac{\cov[\rv{y}_{\mi{M}}]{\rv{y}_{\mi{m}}}{\mi{M}}_{l\x l}}{\dev{\rv{y}_{\mi{m}}}{\mi{m}}_{l} \dev{\rv{y}_{\mi{M}}}{\mi{M}}_{l}}
        = \frac{\dev{\rv{y}_{\mi{m}}}{\mi{m}}_{l}}{\dev{\rv{y}_{\mi{M}}}{\mi{M}}_{l}}
    \end{equation}
    The square of the correlation is a discerning measure of the quality of prediction \cite{Chicco2021} called the coefficient of determination. Which is now clearly identical to the closed Sobol' index
    \begin{equation} \label[definition]{def:GSAI:R2}
        \te[l]{\rv{R}_{\mi{mM}}}^{2} = \frac{\cov{\rv{y}_{\mi{m}}}{\mi{m}}_{l\x l}}{\dev{\rv{y}_{\mi{M}}}{\mi{M}}_{l} \dev{\rv{y}_{\mi{M}}}{\mi{M}}_{l}} \deqr \te[l]{\rv{S}_{\mi{m}}}
    \end{equation}
    Sobol' indices may thus be used to identify reduced models $\tte[l]{\rv{y}_{\mi{m}}}$ which adequately mimic the full model $\tte[l]{\rv{y}_{\mi{M}}}$.
    A closed index close to 1 confirms that the two models make nearly identical predictions. Simplicity and economy (not least of calculation) motivate the adoption of a reduced model, a closed Sobol' index close to 1 is what justifies it. This is precisely equivalent to screening out (obscuring) $\mi{M-m}$ input axes on the grounds that their influence on $\tte[l]{\rv{y}_{\mi{M}}}$ -- measured by their total index $\te[l]{S^{T}_\mi{M-m}}$ -- is close to 0.

\section{Sobol' matrices}\label{sec:GSAM}
    This Section introduces a tensor RV called the Sobol' matrix of an mQU. This is defined by tensorizing the Definitions of the previous Section, replacing $l$ for a scalar output with $\mi{L\x L}$ for a vector output. All division is performed elementwise, inverting Hadamard multiplication $\circ$
    \begin{equation}\label[definition]{def:GSAM:div}
        \te[\square]{*} = \frac{\te[\square]{\bullet}}{\te[\square]{\star}} \quad \Longleftrightarrow
        \quad \te[\square]{*} \circ \te[\square]{\star} = \te[\square]{\bullet} \quad \Longleftrightarrow
        \quad \te[\#]{*} \te[\#]{\star} = \te[\#]{\bullet} \quad \forall \#\in\square
    \end{equation}

    The tensorization of \Cref{def:GSAI:S_m} is problematic, as it may invoke division by zero. Instead, we shall tensorize the alternative \Cref{def:GSAI:R2} of a closed Sobol' index as a coefficient of determination. We shall then compare the two definitions in the light of a toy example.

    The closed Sobol' matrix of vector output $\tte[\mi{L}]{y(\tte[\mi{M+1}]{\rv{u}})}$ with respect to input axes $\mi{m}$ is defined as the tensor RV
    \begin{equation}\label[definition]{def:GSAM:S_m}
            \te[\mi{L\x L}]{\rv{S}_\mi{m}} 
            \deq \frac{\cov{\rv{y}_{\mi{m}}}{\mi{m}}_{\mi{L\x L}}}{\dev{\rv{y}_{\mi{M}}}{\mi{M}}_{\mi{L}} \otimes \dev{\rv{y}_{\mi{M}}}{\mi{M}}_{\mi{L}}}
    \end{equation}
    The complement of a closed Sobol' matrix is called a total Sobol' matrix
    \begin{equation}\label[definition]{def:GSAM:ST_m}
        \te[\mi{L\x L}]{\rv{S}^{T}_\mi{M-m}} \deq \te[\mi{L\x L}]{\rv{S}_{\mi{M}}} - \te[\mi{L\x L}]{\rv{S}_\mi{m}}
    \end{equation}
    Let us emphasise that these Definitions are equivalent to the Sobol' index \Cref{def:GSAI:S_m,def:GSAI:ST_m} on the Sobol' matrix diagonal
    \begin{equation*}
        \te[l\x l]{\rv{S}_\mi{m}} =\te[l]{\rv{S}_\mi{m}} \QT{and} \te[l\x l]{\rv{S}^{T}_\mi{m}} =\te[l]{\rv{S}^{T}_\mi{m}}
        \qquad \forall l\in\mi{L}
    \end{equation*}

    Just as a Sobol' index is a proportion of output variance, each element of a Sobol' matrix is a component of output correlation. Unlike the former, the latter can be negative.
    It is still sometimes useful to examine the first-order matrix for a single input axis $m$
    \begin{equation}\label[definition]{def:GSAM:S_mp}
        \te[\mi{L\x L}]{\rv{S}_{m}} \deq \te[\mi{L\x L}]{\rv{S}_\mi{m+1-m}}
    \end{equation}
    However, there is no ordering analagous to \Cref{eq:GSAI:order} as contributions of opposing sign may offset each other. 

    Closed Sobol' matrices are related to the marginalized correlation matrices defined by tensorizing \Cref{def:GSAI:R_m}
    \begin{equation}\label[definition]{def:GSAM:R_mmp}
        \te[\mi{L\x L^{\prime}}]{\rv{R}_{\mi{mm^{\prime}}}} 
        \deq \frac{\cov[\rv{y}_{\mi{m}^{\prime}}]{\rv{y}_{\mi{m}}}{\mi{m}}_{\mi{L\x L^{\prime}}}}{\dev{\rv{y}_{\mi{m}}}{\mi{m}}_{\mi{L}} \otimes \dev{\rv{y}_{\mi{m^{\prime}}}}{\mi{m^{\prime}}}_{\mi{L^{\prime}}}}
        \qquad \mi{L^{\prime}}\deq \mi{L} \quad \forall \mi{m},\mi{m^{\prime}} \in \mi{M}
    \end{equation}
    In terms of which
    \begin{equation}\label[definition]{def:GSAM:S_R}
            \te[l\x l^{\prime}]{\rv{S}_\mi{m}} = \te[l\x l^{\prime}]{\rv{R}_{\mi{mm}}} 
            \sqrt{\te[l\x l]{\rv{S}_\mi{m}}\te[l^{\prime}\x l^{\prime}]{\rv{S}_\mi{m}}} \qquad \forall l, l^{\prime} \in \mi{L}
    \end{equation}
    
    Consistency demands, for any $l\in{\mi{L}}$ and any realization $\tte[M]{\rv{u}}=\tte[M]{u}$
    \begin{subequations}\label{eq:GSAM:order}
        \begin{gather}
            \te[\mi{L\x L}]{\rv{S}_\mi{0}}=\te[\mi{L\x L}]{\rv{S}^{T}_\mi{0}} = \te[\mi{L\x L}]{0} \\
            \te[\mi{L\x L}]{\rv{S}_\mi{M}}=\te[\mi{L\x L}]{\rv{S}^{T}_\mi{M}} =\te[\mi{L\x L}]{\rv{R}_{\mi{MM}}} \\
            \te[l\x l]{0} \leq
            \modulus{\te[l\x l^{\prime}]{\rv{S}_\mi{m}}} \leq
            \sqrt{\te[l\x l]{\rv{S}_\mi{m}}\te[l^{\prime}\x l^{\prime}]{\rv{S}_\mi{m}}} \leq 1 \qquad \forall l, l^{\prime} \in \mi{L}
        \end{gather}
    \end{subequations}
    So the void model $\tte[\mi{L}]{\rv{y}_{\mi{0}}}$ still explains nothing, the full model $\tte[\mi{L}]{\rv{y}_{\mi{M}}}$ still explains everything explicable, and every other mQU lies somewhere in between. The last inequality shows that inputs can only influence the correlation between two outputs they are relevant to.

    Let us pursue an illuminating toy example to illustrate these definitions. Take some even $M$ and define
    \begin{equation*}
            w \colon \uniti^{M/2} \to \st{R} \QT{such that} \ev{w}{\mi{M/2}} = 0 \QT{and} \cov{w}{\mi{M/2}} = 1
    \end{equation*}
    to construct the ($L=2$) MQU
    \begin{equation*}
        y(\te[\mi{M+1}]{\rv{u}}) = 
        \begin{bmatrix}
            w(\te[\mi{M/2}]{\rv{u}}) + w(\te[\mi{M-M/2}]{\rv{u}}) \\
            w(\te[\mi{M/2}]{\rv{u}}) - w(\te[\mi{M-M/2}]{\rv{u}})
        \end{bmatrix}
    \end{equation*}
    The Sobol' matrix of the full model is easily calculated
    \begin{equation*}
        {\rv{S}_\mi{M}} = 
        \begin{bmatrix}
            1 & 0 \\
            0 & 1
        \end{bmatrix}
    \end{equation*}
    Two important reduced models have the Sobol' matrices
    \begin{equation*}
        {\rv{S}_\mi{M/2}} = {\rv{S}^{T}_\mi{M/2}} = 
        \begin{bmatrix}
            1/2 & 1/2 \\
            1/2 & 1/2
        \end{bmatrix}
        \QT{;}
        {\rv{S}_\mi{M-M/2}} = {\rv{S}^{T}_\mi{M-M/2}} = 
        \begin{bmatrix}
            \phantom{-}1/2 & -1/2 \\
            -1/2 & \phantom{-}1/2
        \end{bmatrix}
    \end{equation*}
    The Sobol' matrices reveal that input axes $\mi{M/2}$ influence both outputs in the same sense, while input axes $\mi{M-M/2}$ influence the two outputs in opposite senses.
    In real-world examples, this information is often valuable. If the two components of $y$ are the yield and the purity of a pharamaceutical product one would initially focus on inputs $\mi{M/2}$ when trying to maximise both outputs simultaneously. Alternatively, if the two components of $y$ are the efficacy and the cost of an industrial process, one would initally vary inputs $\mi{M-M/2}$ in an effort to improve the tradeoff between the two outputs.
    Such prescriptions are particularly simple in this case because the total Sobol' matrices are equal to the closed Sobol' matrices, revealing that inputs $\mi{M/2}$ do not cooperate with inputs $\mi{M-M/2}$.

    The direct tensorization of Sobol' index \Cref{def:GSAI:S_m} would define the matrix
    \begin{equation*}
            \te[\mi{L\x L}]{\rv{\hat{S}}_\mi{m}} \deq \frac{\te[\mi{L\x L}]{\rv{S}_\mi{m}}}{\te[\mi{L\x L}]{\rv{S}_\mi{M}}} \QT{which is undefined whenever} \te[l\x l^{\prime}]{\rv{S}_\mi{M}}=0 \T{ for some } l,l^{\prime} \in \mi{L}
    \end{equation*}
    In our toy example, the most interesting information -- the off-diagonal indices for two important reduced models -- is undefined and completely lost by this definition.

    We close this Section by considering how a Sobol' matrix may be summarized in a single number. This is often desired to present ``the'' sensitivity of multi-output to a reduced set of inputs $\mi{m}$. The simple answer is to define a seminorm on the Sobol' matrix \cite[pp.314]{Schechter1997}. The seminorm should be chosen according to one's interest in the multi-output. For example, one could use the determinant $\modulus{\tte[\mi{l\x l}]{\rv{S}_{\mi{m}}}}$ of a chosen submatrix $\mi{l}\subseteq\mi{L}$, which is just the modulus in case $\mi{l}$ is singleton. The Sobol' matrix provides a platform for investigating a variety of such measures. However, all matrices and measures remain RVs at this point, and the next two Sections are devoted to extracting statistics -- i.e. determined quantities -- from these.

\section{Sobol matrix statistics}\label{sec:GSAS}
    At this stage, Sobol' matrices have been robustly defined as tensor RVs -- i.e. functions of the undetermined input $\tte[M]{\rv{u}}$. This Section extracts two determined matrix statistics from each RV: its expected value and its standard deviation or error.

    Expected values are written as the italic version of the underlying RV, starting with
    \begin{subequations}\label[definition]{def:GSAS:VD}
        \begin{align}
            \QT{marginal variance} \te[\mi{L\x L}]{V_\mi{m}} &\deq \evt{\,\cov{\rv{y}_{\mi{m}}}{\mi{m}}}{M}_{\mi{L\x L}}  \label[definition]{def:GSAS:V} \\
            \QT{marginal deviation} \te[\mi{L}]{D_\mi{m}} &\deq \evt{\,\dev{\rv{y}_{\mi{m}}}{\mi{m}}}{M}_{\mi{L}}
            \quad \Longrightarrow \quad
            \te[l]{D_\mi{m}} = \sqrt{\te[l\x l]{V_\mi{m}}} \quad \forall l\in\mi{L} \label[definition]{def:GSAS:D}
        \end{align}
    \end{subequations}
    The GSA statistics used to assess input relevance are the expected Sobol' matrices
    \begin{subequations}\label[definition]{def:GSAS:SST}
        \begin{align}
            \QT{closed Sobol' matrix}\te[\mi{L}^{2}]{S_\mi{m}} &\deq \evt{\te[\mi{L}^{2}]{\rv{S}_\mi{m}}}{M} 
            = \frac{\te[\mi{L}^{2}]{V_\mi{m}}}{\te[\mi{L}]{D_\mi{M}} \otimes \te[\mi{L}]{D_\mi{M}}} \\    
            \QT{total Sobol' matrix}\te[\mi{L}^{2}]{S^{T}_\mi{M-m}} &\deq \evt{\te[\mi{L}^{2}]{\rv{S}^{T}_\mi{M-m}}}{M} = \te[\mi{L}^{2}]{S_\mi{M}} - \te[\mi{L}^{2}]{S_\mi{m}}
        \end{align}
    \end{subequations}

    In order to express standard errors -- due to undetermined input $\tte[M]{\rv{u}}$ -- the notation introduced in \Cref{def:GSAI:D_m} is extended to define the matrix standard deviation $\devt{\te[\mi{L\x L}]{\bullet}}{M}$ of a matrix RV $\tte[\mi{L\x L}]{\bullet}$
    \begin{equation}\label[definition]{def:GSAS:D_m}
        \devt{\te[l\x l^{\prime}]{\bullet}}{M} \deq \sqrt{\cov{\te[l\x l^{\prime}]{\bullet}}{M}} \qquad \forall l,l^{\prime} \in \mi{L}
    \end{equation}
    The GSA statistics used to assess the standard error of Sobol' matrices are
    \begin{subequations}\label[definition]{def:GSAS:TTT}
        \begin{align}
            \QT{closed Sobol' matrix error}\te[\mi{L}^{2}]{T_\mi{m}} &\deq \devt{\te[\mi{L}^{2}]{\rv{S}_\mi{m}}}{M} \label[definition]{def:GSAS:T}\\
            \QT{total Sobol' matrix error}
            \te[\mi{L}^{2}]{T_{\mi{M-m}}^{T}} &\deq 
            \te[\mi{L}^{2}]{T_{\mi{M}}} + \te[\mi{L}^{2}]{T_{\mi{m}}}
            \geq \devt{\,\te[\mi{L}^{2}]{\rv{S}_{\mi{M-m}}^{T}}}{M}
        \end{align}
    \end{subequations}
    The total Sobol' matrix error is a conservative statistic which achieves equality on the diagonal, and is robust and sufficiently precise for most practical purposes off the diagonal.
    Attempts at greater precision are apt to yield nonsense when implemented numerically.
    This is because the total Sobol' matrix \Cref{def:GSAM:ST_m} is the difference between two terms which are often highly correlated with each other. Tiny differences between correlated terms are swamped by numerical error in implementation, so the resulting computation is wildly unreliable.

    To calculate the closed Sobol' matrix error $\tte[\mi{L\x L}]{T_\mi{m}}$ we use the Taylor series method \cite[pp.353]{Kendall1994}, which is valid provided $\covt{\,\te[l\x l]{\rv{y}_{\mi{M}}}}{\mi{M}}$ is well approximated by its mean
    \begin{equation*}
        \te[l\x l]{V_{\mi{M}}} \gg \big\vert\covt{\,\te[l\x l]{\rv{y}_{\mi{M}}}}{\mi{M}} - \te[l\x l]{V_{\mi{M}}}\big\vert
    \end{equation*}
    This is essentially the positive variance clause in the \Cref{def:MQU:y} of an MQU, designed to prohibit constant or pure noise outputs which do not depend on any of the determined input axes $\mi{M}$.

    The covariance between two marginal covariance matrix RVs is the determined 4-tensor
    \begin{equation}\label[definition]{def:GSAS:W}
        \te[\mi{L}^{4}]{W_{\mi{mm^{\prime}}}} \deq \cov[\cov{\rv{y}_{\mi{m^{\prime}}}}{\mi{m^{\prime}}}_{\mi{L\x L}}]{\cov{\rv{y}_{\mi{m}}}{\mi{m}}_{\mi{L\x L}}}{M} \qquad \forall \mi{m},\mi{m^{\prime}} \in \mi{M}
    \end{equation}
    used to define the determined matrix $\tte[\mi{L\x L}]{Q_{\mi{m}}}$ elementwise
    \begin{equation}\label[definition]{def:GSAS:Q}
        \begin{aligned}
            \te[l\x l^{\prime}]{Q_{\mi{m}}} &\deq 
            \te[(l\x l^{\prime})^{2}]{W_{\mi{mm}}} \\
            &\phantom{\deq}- \te[l\x l^{\prime}]{V_{\mi{m}}} \sum_{l^{\ddag} \in \set{l, l^{\prime}}} \frac{\te[l^{\ddag}\x l^{\ddag}\x l\x l^{\prime}]{W_{\mi{mM}}}} {\te[l^{\ddag}]{D_{\mi{M}}}^{2}} \\                
            &\phantom{\deq}+ \frac{\te[l\x l^{\prime}]{V_{\mi{m}}}^{2}}{4} \sum_{l^{\ddag} \in \set{l, l^{\prime}}} 
            \frac{\te[l^{\ddag}\x l^{\ddag}\x l\x l]{W_{\mi{MM}}}} {\te[l^{\ddag}]{D_{\mi{M}}}^{2} \te[l]{D_{\mi{M}}}^{2}}
            + \frac{\te[l^{\ddag}\x l^{\ddag}\x l^{\prime}\x l^{\prime}]{W_{\mi{MM}}}} {\te[l^{\ddag}]{D_{\mi{M}}}^{2} \otimes \te[l^{\prime}]{D_{\mi{M}}}^{2}} \\
        \end{aligned}
    \end{equation}
    The Taylor series estimate of the closed Sobol matrix error \Cref{def:GSAS:T} is finally calculated elementwise as
    \begin{equation} \label{eq:GSAS:Tm}
        \te[l\x l^{\prime}]{T_{\mi{m}}} \deq 
        \devt{\,\te[l\x l^{\prime}]{\rv{S}_{\mi{m}}}}{M} = \frac{\sqrt{\te[l\x l^{\prime}]{Q_{\mi{m}}}}}{\te[l]{D_{\mi{M}}} \te[\mi{l^{\prime}}]{D_{\mi{M}}}}
    \end{equation}
    It is satisfying to note that these equations enforce
    \begin{equation*} \label{eq:GSAS:Tm:diag}
        \te[l\x l]{T_{\mi{M}}}=0 \QT{on the diagonal} \te[l\x l]{\rv{S}_{\mi{M}}}=1        
    \end{equation*}

\section{Sobol' matrix statistics from MQU moments}\label{sec:Mom}
    The central problem in calculating errors on Sobol' matrices is that they involve ineluctable covariances over undetermined noise between differently marginalized SPs. But marginalization and covariance determination are both a matter of taking expectations. Using Fubini's theorem \cite[pp.77]{Williams1991} the ineluctable can be avoided by reversing the order of expectations -- taking moments over ungoverned noise, then marginalizing.

    \Cref{sec:MQU} invoked the Kolmogorov extension theorem \cite[pp.124]{Rogers.Williams2000} to guarantee the existence of a tensor-valued SP called a finite-dimensional distribution 
    \begin{equation*}
        \te[\mi{L\x N}]{y(\te[\mi{M+1\x N}]{\rv{u}})} \QT{for any} N\in\st{Z}^{+}
    \end{equation*}
    This used to construct MQU moments as follows. The first moment is the determined statistic
    \begin{equation}\label{eq:Mom:muM}
        \te[\mi{L}]{\mu_{\mi{M}}} \deq \evt{\te[\mi{L}]{y(\te[\mi{M+1}]{\rv{u}})}}{M}                
    \end{equation}
    engendering the centralized MQU
    \begin{equation}\label[definition]{def:Mom:e}
        \te[\mi{L}]{\rv{e}_{\mi{M}}} \deq \te[\mi{L}]{y(\te[\mi{M+1}]{\rv{u}})} - \te[\mi{L}]{\mu_{\mi{M}}}
    \end{equation}
    The $N^{\mathrm{th}}$ MQU moment is the determined statistic
    \begin{equation}\label{eq:Mom:muMMM}
        \te[\mi{L}^{N}]{\mu_{\mi{M\ldots M}^{n\prime}\mi{\ldots M}^{(N-1)\prime}}} \deq \ev{\tte[\mi{L}]{\rv{e}_{\mi{M}}}\otimes\cdots\otimes\tte[\mi{L}^{n\prime}]{\rv{e}_{\mi{M}^{n\prime}}}\otimes
        \cdots\otimes\tte[\mi{L}^{(N-1)\prime}]{\rv{e}_{\mi{M}^{(N-1)\prime}}}}{M}
    \end{equation}
    where, to be clear
    \begin{equation*}
        \begin{aligned}
            n &\in \mi{N} \\
            \mi{m}^{n\prime} \subseteq \mi{M}^{n\prime} &\deq \mi{M} \supseteq \mi{m} \\
            \mi{l}^{n\prime} \subseteq \mi{L}^{n\prime} &\deq \mi{L} \supseteq \mi{l}            
        \end{aligned}
    \end{equation*}
    Throughout this work, primes on uppercase constants are for bookkeeping only, they do not affect the value. 
    Primes on lowercase variables do affect the value, so it may or may not be the case that $\mi{m}^{n\prime}=\mi{m}$ for example. An $N^{\mathrm{th}}$ mQU moment is a determined statistic obtained by marginalizing an $N^{\mathrm{th}}$ MQU moment
    \begin{equation}\label{eq:Mom:mummm}
        \begin{aligned}
            \te[\mi{L}^{N}]{\mu_{\mi{m}\mi{\ldots m}^{(N-1)\prime}}} &\deq \evt{\ldots\evt{\,\ev{\tte[\mi{L}]{\rv{e}_{\mi{M}}}\otimes\cdots\otimes\tte[\mi{L}^{(N-1)\prime}]{\rv{e}_{\mi{M}^{(N-1)\prime}}}}{M}}{\mi{M}^{(N-1)\prime}-\mi{m}^{(N-1)\prime}}}{\mi{M-m}} \\
            &\deqr \ev{\tte[\mi{L}]{\rv{e}_{\mi{m}}}\otimes\cdots\otimes\tte[\mi{L}^{(N-1)\prime}]{\rv{e}_{\mi{m}^{(N-1)\prime}}}}{M}            
        \end{aligned}
    \end{equation}
    The last line relies on Fubini's theorem. The law of iterated expectations further entails
    \begin{equation}\label{eq:Mom:reduction}
        \te[\mi{L}^{N}]{\mu_{\mi{0\ldots 0}\mi{m}^{n\prime}\mi{\ldots m}^{(N-1)\prime}}} 
        = \evt{\te[\mi{L}^{N}]{\mu_{\mi{m\ldots m}\mi{m}^{n\prime}\mi{\ldots m}^{(N-1)\prime}}}}{\mi{m}} 
        = \evt{\te[\mi{L}^{N}]{\mu_{\mi{M\ldots M}\mi{m}^{n\prime}\mi{\ldots m}^{(N-1)\prime}}}}{\mi{M}}
    \end{equation}
    This reduction will be used repeatedly in the remainder of this Section. Although $\te[\mi{L}^{2}]{\mu_{\mi{00}}} \geq 0$ in general, \Cref{def:MQU:y,def:GSAI:V_m,def:GSAS:VD} conveniently zero the fully marginalized quantities
    \begin{equation} \label{eq:Mom:V0}
        \covt{\te[\mi{L}^{2}]{\rv{y}_{\mi{0}}}}{\mi{0}} = \te[\mi{L}^{2}]{V_{\mi{0}}} = \te[\mi{L}]{\mu_{\mi{0}}}^{2} = \te[\mi{L}^2]{0}
    \end{equation}
    
    The expected conditional variance in \Cref{def:GSAS:VD} amounts to
    \begin{equation}\label{eq:Mom:V}
        \begin{aligned}
            \te[\mi{L}^{2}]{V_{\mi{m}}} &\deq \evt{\covt{\te[\mi{L}^{2}]{\rv{y}_{\mi{m}}}}{\mi{m}}}{M} 
            = \evt{\ev{\te[\mi{L}]{\rv{y}_{\mi{m}}}^{2} - \te[\mi{L}]{\rv{y}_{\mi{0}}}^{2}}{\mi{m}}}{M} \\
            &\phantom{:}= \evt{\ev{\te[\mi{L}]{\rv{e}_{\mi{m}} + \mu_{\mi{m}}}^{2}}{M}}{\mi{m}}
            - \ev{\te[\mi{L}]{\rv{e}_{\mi{0}} + \mu_{\mi{0}}}^{2}}{M} \\
            &\phantom{:}= \ev{\te[\mi{L}^2]{\mu_{\mi{mm}}} + \te[\mi{L}]{\mu_{\mi{m}}}^{2}}{\mi{m}}  - \te[\mi{L}^2]{\mu_{\mi{00}}} - \te[\mi{L}]{\mu_{\mi{0}}}^{2} \\
            &\phantom{:}= \ev{\te[\mi{L}]{\mu_{\mi{m}}}^{2}}{\mi{m}}
        \end{aligned}
    \end{equation}
    Substituting this formula in \Cref{def:GSAS:V,def:GSAS:SST} determines the closed and total Sobol' matrices in terms of mQU moments.

    The covariance between conditional variances in \Cref{def:GSAS:W} is
    \begin{equation*}
        \begin{aligned}
            \te[\mi{L}^{2}\x \mi{L^{\prime}}^{2}]{W_{\mi{mm^{\prime}}}} &\deq \cov[\covt{\te[\mi{L^{\prime}}^{2}]{\rv{y}_{\mi{m^{\prime}}}}}{\mi{m^{\prime}}}]{\covt{\te[\mi{L}^{2}]{\rv{y}_{\mi{m}}}}{\mi{m}}}{M} \\
            &\phantom{:}=
            \cov[\ev{\te[\mi{L^{\prime}}]{\rv{y}_{\mi{m^{\prime}}}}^{2} - \te[\mi{L^{\prime}}]{\rv{y}_{\mi{0}}}^{2}}{\mi{m^{\prime}}}]{\ev{\te[\mi{L}]{\rv{y}_{\mi{m}}}^{2} - \te[\mi{L}]{\rv{y}_{\mi{0}}}^{2}}{\mi{m}}}{M} \\
            &\phantom{:}=
            \ev{\ev{\te[\mi{L}]{\rv{y}_{\mi{m}}}^{2} - \te[\mi{L}]{\rv{y}_{\mi{0}}}^{2}}{\mi{m}} \otimes\ev{\te[\mi{L^{\prime}}]{\rv{y}_{\mi{m^{\prime}}}}^{2} - \te[\mi{L^{\prime}}]{\rv{y}_{\mi{0}}}^{2}}{\mi{m^{\prime}}}}{M} - \te[\mi{L}^2]{V_{\mi{m}}}\otimes \te[\mi{L^{\prime}}^2]{V_{\mi{m^{\prime}}}} \\       
            &\phantom{:}= \te[\mi{L}^2\x \mi{L^\prime}^2]{A_{\mi{mm^{\prime}}}-A_{\mi{0m^{\prime}}}-A_{\mi{m0}}+A_{\mi{00}}}
        \end{aligned}
    \end{equation*}
    where
    \begin{equation*}
        \begin{aligned}
            \te[\mi{L}^2\x \mi{L^\prime}^2]{A_{\mi{mm^{\prime}}}}
            &\deq \evt{\evt{\ev{\te[\mi{L}]{\rv{y}_{\mi{m}}}^{2} \otimes \te[\mi{L^{\prime}}]{\rv{y}_{\mi{m^{\prime}}}}^{2}}{\mi{m}}}{\mi{m^{\prime}}}}{M} - \te[\mi{L}^2]{V_{\mi{m}}}\otimes \te[\mi{L^{\prime}}^2]{V_{\mi{m^{\prime}}}} \\
            &\phantom{:}= \evt{\evt{\ev{
                \te[\mi{L}]{\rv{e}_{\mi{m}}+\mu_{\mi{m}}}^{2} \otimes \te[\mi{L^{\prime}}]{\rv{e}_{\mi{m^{\prime}}}+ \mu_{\mi{m^{\prime}}}}^{2} - 
                \te[\mi{L}]{\mu_{\mi{m}}}^{2} \otimes \te[\mi{L^{\prime}}]{\mu_{\mi{m^{\prime}}}}^{2}
            }{M}}{\mi{m^{\prime}}}}{\mi{m}}
        \end{aligned}
    \end{equation*}
    exploiting the fact that $V_{\mi{0}} = \te[\mi{L}^2]{0}$. \Cref{eq:Mom:reduction} cancels all terms beginning with $\te[\mi{L}]{\rv{e}_{\mi{m}}}^{2}$, across $A_{\mi{mm^{\prime}}}-A_{\mi{0m^{\prime}}}-A_{\mi{m0}}+A_{\mi{00}}$. All remaining terms ending in $\te[\mi{L^{\prime}}]{\mu_{\mi{m^{\prime}}}}^{2}$ are eliminated by centralization $\evt{\,\tte[]{\rv{e}_{\mi{m}}}}{M} = 0$.
    Similar arguments eliminate $\te[\mi{L^{\prime}}]{\rv{e}_{\mi{m^{\prime}}}}^{2}$ and $\te[\mi{L}]{\mu_{\mi{m}}}^{2}$.
    Effectively then
    \begin{equation*}
        \te[\mi{L}^4]{A_{\mi{mm^{\prime}}}} = \sum_{\pi(\mi{L}^{2})} \sum_{\pi(\mi{L^{\prime}}^{2})}
        \evt{\evt{\te[\mi{L}^{2} \x \mi{L^{\prime}}^{2}]{\mu_{\mi{m}} \otimes \mu_{\mi{mm^{\prime}}} \otimes \mu_{\mi{m^{\prime}}}}}{\mi{m^{\prime}}}}{\mi{m}}
    \end{equation*}
    so \Cref{eq:Mom:V0} entails
    \begin{equation}\label{eq:Mom:W}
        \te[\mi{L}^4]{W_{\mi{mm^{\prime}}}} = \sum_{\pi(\mi{L}^{2})} \sum_{\pi(\mi{L^{\prime}}^{2})}
        \evt{\evt{\te[\mi{L}^{2} \x \mi{L^{\prime}}^{2}]{\mu_{\mi{m}} \otimes \mu_{\mi{mm^{\prime}}} \otimes \mu_{\mi{m^{\prime}}}}}{\mi{m^{\prime}}}}{\mi{m}}
    \end{equation}
    where each summation is over permutations of tensor axes
    \begin{equation*}
        \pi(\mi{L}^{2}) \deq \set{(\mi{L}\x\mi{L^{\prime\prime}}), (\mi{L^{\prime\prime}}\x\mi{L})} \QT{;} \pi(\mi{L^{\prime}}^{2}) \deq \set{(\mi{L^{\prime}}\x\mi{L^{\prime\prime\prime}}), (\mi{L^{\prime\prime\prime}}\x\mi{L^{\prime}})}
    \end{equation*}
    Substituting \Cref{eq:Mom:W} in \Cref{def:GSAS:Q,def:GSAS:TTT} determines the closed and total Sobol' matrix errors in terms of marginalized MQU moments.

    \Cref{eq:Mom:W} has an important property. Any constant term $E$ added to the MQU variance $\mu_{\mi{mm^{\prime}}}\mapsto \mu_{\mi{mm^{\prime}}}+E$ adds nothing to $W_\mi{mm^{\prime}}$ by reduction \Cref{eq:Mom:reduction}. Adding homoskedastic noise to an MQU directly affects neither $S_\mi{m}$ nor $T_\mi{m}$. GSA compares a reduced model to a full model, largely unaffected by noise common to both.
    
    So Sobol' matrices may be calculated from the first MQU moment alone. Their standard errors may be calculated from the first two MQU moments alone.
    The first MQU moment $\tte[\mi{L}]{\mu_{\mi{M}}}$ is simply the mean prediction of the MQU, a function of $\tte[\mi{M}]{u}$. The second MQU moment $\tte[\mi{L}^{2}]{\mu_{\mi{MM}}}$ is known as the covariance function or kernel of an SP
    \begin{equation} \label{eq:Mom:k}
        \te[\mi{L\x L^{\prime}}]{\mu_{\mi{MM^{\prime}}}} = \te[\mi{L\x L^{\prime}}]{k(\tte[\mi{M}]{u},\tte[\mi{M^{\prime}}]{u^{\prime}})} \deq 
    \cov[{\te[\mi{L^{\prime}}]{y(\te[\mi{M^{\prime}+1}]{\rv{u}^{\prime}})}}]{\te[\mi{L}]{y(\te[\mi{M+1}]{\rv{u}})}}{M}
    \end{equation}
    It expresses the linkage between MQU predictions at two determined inputs, and is central to the theory and implementation of GPs \cite{Rasmussen.Williams2005}.

    Analytic expressions for Sobol' matrices and their standard errors are given in \cite{MiltonRobert.etal.2025} for GPs with a Radial Basis Function (RBF) kernel \cite{Rasmussen.Williams2005}. These have been implemented in Python as the RomCom package \cite{RomCom2025}, which is used in the next two Sections to benchmark the calculation of Sobol' matrices and their standard errors from MQU moments.

\section{Test functions}\label{sec:Func}
    The SALib package \cite{SALib2025} implements a number of scalar test functions for benchmarking sensitivity analysis. This study uses 3 of these
    \begin{subequations}
        \begin{align}
            \QT{Ishigami \cite{Ishigami}:} &\ish(\tte[\mi{3}]{u};A,B) \deq \left(1 + B\tte[2]{u}^{4} \right)\sin(\tte[0]{u}) + A \sin(\tte[1]{u})^{2} \\
            \QT{Sobol' G \cite{Saltelli2010}:} &\sob(\tte[\mi{5}]{u};\tte[\mi{5}]{A},\tte[\mi{5}]{B}) \deq \prod_{m=0}^{4} \frac{(1+\tte[m]{B})\modulus{2\tte[m]{u}-1}^{\tte[m]{B}}+\tte[m]{A}}{1 + \tte[m]{A}}\\
            \QT{Oakley 2004 \cite{Oakley.OHagan2004}:} &\oak(\tte[\mi{5}]{u};\tte[\mi{5}]{A},\tte[\mi{5\x 5}]{B}) \deq \tte[\mi{5}]{u} \tte[\mi{5}]{A} + \tte[\mi{5}]{u} \tte[\mi{5\x 5^{\prime}}]{B} \tte[\mi{5^{\prime}}]{u}
        \end{align}
    \end{subequations}
    The last of these uses the Einstein summation convention for tensor contraction. It is just a quadratic form in $\tte[\mi{5}]{u}$ parametrized by the matrix $\tte[\mi{5\x 5}]{B}$ plus a linear term parametrized by the vector $\tte[\mi{5}]{A}$.
    To parametrize the (modified) Sobol' G function and (restricted) Oakley 2004 function we define the vectors
    \begin{equation}
        \tte[\mi{5}]{A_{g}} \deq \left\lbrack
            \begin{array}{r}
            3 \\
            6 \\
            9 \\
            18 \\
            27 \\
            \end{array}
            \right\rbrack
        \QT{;}
        \tte[\mi{5}]{A_{G}} \deq \left\lbrack
            \begin{array}{r}
             {1}/{2} \\
             1 \\
             2 \\
             4 \\
             8 \\
            \end{array}
            \right\rbrack
        \QT{;}
        \tte[\mi{5}]{A_{O}} \deq \left\lbrack
            \begin{array}{r}
                5 \\
                {35}/{8} \\
                {15}/{4} \\
                {25}/{8} \\
                {5}/{2} \\
            \end{array}
            \right\rbrack
        \end{equation}
    and the matrices
    \begin{equation}
            \tte[\mi{5\x 5}]{B_{+}} \deq \left\lbrack
                \begin{array}{rrrrr}
                 5 & {29}/{6} & {14}/{3} & {9}/{2} & {13}/{3} \\
                 {25}/{6} & 4 & {23}/{6} & {11}/{3} & {7}/{2} \\
                 {10}/{3} & {19}/{6} & 3 & {17}/{6} & {8}/{3} \\
                 {5}/{2} & {7}/{3} & {13}/{6} & 2 & {11}/{6} \\
                 {5}/{3} & {3}/{2} & {4}/{3} & {7}/{6} & 1 \\
                \end{array}
                \right\rbrack
            \QT{;}
            \tte[m \x m^{\prime}]{B_{-}} \deq \tte[4-m \x 4-m^{\prime}]{B_{+}}
    \end{equation}
    The $L=9$ MNU used for benchmarking is
    \begin{equation}
        f(\tte[\mi{M}]{u}) \deq \left\lbrack
            \begin{array}{l}
                \ish(2 \pi \tte[\mi{3}]{u} - \pi \tte[\mi{3}]{1}; 7, 1/10) \\
                \ish(2 \pi \tte[\mi{3}]{u} - \pi \tte[\mi{3}]{1}; 20, 1) \\
                \ish(2 \pi \tte[\mi{3}]{u} - \pi \tte[\mi{3}]{1}; 0, 0) \\
                \sob(\tte[\mi{5}]{u}; \tte[\mi{5}]{A_{g}}, \tte[\mi{5}]{2}) \\
                \sob(\tte[\mi{5}]{u}; \tte[\mi{5}]{A_{G}}, \tte[\mi{5}]{2}) \\
                \sob(\tte[\mi{5}]{u}; \tte[\mi{5}]{A_{G}}, \tte[\mi{5}]{4}) \\
                \oak(2\tte[\mi{5}]{u} - \tte[\mi{5}]{1}; \tte[\mi{5}]{A_{O}}, \tte[\mi{5\x 5}]{0}) \\
                \oak(2\tte[\mi{5}]{u} - \tte[\mi{5}]{1}; \tte[\mi{5}]{A_{O}}, \tte[\mi{5\x 5}]{B_{+}}) \\
                \oak(2\tte[\mi{5}]{u} - \tte[\mi{5}]{1}; -\tte[\mi{5}]{A_{O}}, \tte[\mi{5\x 5}]{B_{-}})
            \end{array}
            \right\rbrack
    \end{equation}
    The correlation matrix of this MNU is
    \begin{equation*}
        {S_{\mi{5}}} = \left\lbrack\begin{array}{rrrrrrrrr}
            1.000 & 0.896 & 0.560 & -0.073 & -0.078 & -0.131 & 0.254 & 0.125 & -0.159 \\
            0.896 & 1.000 & 0.593 & -0.032 & -0.034 & -0.057 & 0.268 & 0.146 & -0.161 \\
            0.560 & 0.593 & 1.000 & 0.000 & 0.000 & 0.000 & 0.453 & 0.264 & -0.264 \\
            -0.073 & -0.032 & 0.000 & 1.000 & 0.944 & 0.825 & 0.000 & 0.251 & 0.136 \\
            -0.078 & -0.034 & 0.000 & 0.944 & 1.000 & 0.926 & 0.000 & 0.232 & 0.137 \\
            -0.131 & -0.057 & 0.000 & 0.825 & 0.926 & 1.000 & 0.000 & 0.197 & 0.116 \\
            0.254 & 0.268 & 0.453 & 0.000 & 0.000 & 0.000 & 1.000 & 0.582 & -0.582 \\
            0.125 & 0.146 & 0.264 & 0.251 & 0.232 & 0.197 & 0.582 & 1.000 & 0.206 \\
            -0.159 & -0.161 & -0.264 & 0.136 & 0.137 & 0.116 & -0.582 & 0.206 & 1.000 \\
        \end{array}\right\rbrack
    \end{equation*}
    The Sobol' G functions $\tte[3]{f},\tte[4]{f},\tte[5]{f}$ are invariant under the transformation $\tte[m]{u} \mapsto \tte[m]{1-u}$ -- they possess the reflection symmetries of an even function. The sine function $\tte[2]{f}$ and the linear function $\tte[6]{f}$ preserve magnitude and change sign under the transformation $\tte[m]{u} \mapsto \tte[m]{1-u}$ -- they possess rotational symmetries of an odd function. As a result, the correlation matrix contains a central rectangle of three Sobol' G functions bordered by zero correlations.

    The remaining closed Sobol' matrices of the test MNU are as follows
    \begin{equation*}
        {S_{\mi{4}}} = \left\lbrack\begin{array}{rrrrrrrrr}
            1.000 & 0.896 & 0.560 & -0.073 & -0.078 & -0.131 & 0.254 & 0.125 & -0.159 \\
                0.896 & 1.000 & 0.593 & -0.032 & -0.034 & -0.057 & 0.268 & 0.146 & -0.161 \\
                0.560 & 0.593 & 1.000 & 0.000 & 0.000 & 0.000 & 0.453 & 0.264 & -0.264 \\
                -0.073 & -0.032 & 0.000 & 0.986 & 0.929 & 0.811 & 0.000 & 0.247 & 0.116 \\
                -0.078 & -0.034 & 0.000 & 0.929 & 0.979 & 0.904 & 0.000 & 0.229 & 0.118 \\
                -0.131 & -0.057 & 0.000 & 0.811 & 0.904 & 0.970 & 0.000 & 0.194 & 0.100 \\
                0.254 & 0.268 & 0.453 & 0.000 & 0.000 & 0.000 & 0.916 & 0.533 & -0.533 \\
                0.125 & 0.146 & 0.264 & 0.247 & 0.229 & 0.194 & 0.533 & 0.839 & 0.031 \\
                -0.159 & -0.161 & -0.264 & 0.116 & 0.118 & 0.100 & -0.533 & 0.031 & 0.591 \\
        \end{array}\right\rbrack
    \end{equation*}
    \begin{equation*}
        {S_{\mi{3}}} = \left\lbrack\begin{array}{rrrrrrrrr}
            1.000 & 0.896 & 0.560 & -0.073 & -0.078 & -0.131 & 0.254 & 0.125 & -0.159 \\
                0.896 & 1.000 & 0.593 & -0.032 & -0.034 & -0.057 & 0.268 & 0.146 & -0.161 \\
                0.560 & 0.593 & 1.000 & 0.000 & 0.000 & 0.000 & 0.453 & 0.264 & -0.264 \\
                -0.073 & -0.032 & 0.000 & 0.956 & 0.889 & 0.774 & 0.000 & 0.235 & 0.093 \\
                -0.078 & -0.034 & 0.000 & 0.889 & 0.912 & 0.833 & 0.000 & 0.215 & 0.091 \\
                -0.131 & -0.057 & 0.000 & 0.774 & 0.833 & 0.877 & 0.000 & 0.183 & 0.077 \\
                0.254 & 0.268 & 0.453 & 0.000 & 0.000 & 0.000 & 0.784 & 0.457 & -0.457 \\
                0.125 & 0.146 & 0.264 & 0.235 & 0.215 & 0.183 & 0.457 & 0.622 & -0.096 \\
                -0.159 & -0.161 & -0.264 & 0.093 & 0.091 & 0.077 & -0.457 & -0.096 & 0.359 \\
        \end{array}\right\rbrack
    \end{equation*}
    \begin{equation*}
        {S_{\mi{2}}} = \left\lbrack\begin{array}{rrrrrrrrr}
            0.756 & 0.525 & 0.560 & -0.073 & -0.078 & -0.131 & 0.254 & 0.125 & -0.159 \\
                0.525 & 0.435 & 0.593 & -0.032 & -0.034 & -0.057 & 0.268 & 0.146 & -0.161 \\
                0.560 & 0.593 & 1.000 & 0.000 & 0.000 & 0.000 & 0.453 & 0.264 & -0.264 \\
                -0.073 & -0.032 & 0.000 & 0.848 & 0.765 & 0.661 & 0.000 & 0.202 & 0.059 \\
                -0.078 & -0.034 & 0.000 & 0.765 & 0.741 & 0.660 & 0.000 & 0.181 & 0.057 \\
                -0.131 & -0.057 & 0.000 & 0.661 & 0.660 & 0.664 & 0.000 & 0.154 & 0.048 \\
                0.254 & 0.268 & 0.453 & 0.000 & 0.000 & 0.000 & 0.595 & 0.346 & -0.346 \\
                0.125 & 0.146 & 0.264 & 0.202 & 0.181 & 0.154 & 0.346 & 0.375 & -0.145 \\
                -0.159 & -0.161 & -0.264 & 0.059 & 0.057 & 0.048 & -0.346 & -0.145 & 0.221 \\
        \end{array}\right\rbrack
    \end{equation*}
    \begin{equation*}
        {S_{\mi{1}}} = \left\lbrack\begin{array}{rrrrrrrrr}
            0.314 & 0.332 & 0.560 & 0.000 & 0.000 & 0.000 & 0.254 & 0.148 & -0.148 \\
                0.332 & 0.351 & 0.593 & 0.000 & 0.000 & 0.000 & 0.268 & 0.156 & -0.156 \\
                0.560 & 0.593 & 1.000 & 0.000 & 0.000 & 0.000 & 0.453 & 0.264 & -0.264 \\
                0.000 & 0.000 & 0.000 & 0.632 & 0.515 & 0.438 & 0.000 & 0.139 & 0.028 \\
                0.000 & 0.000 & 0.000 & 0.515 & 0.420 & 0.357 & 0.000 & 0.113 & 0.023 \\
                0.000 & 0.000 & 0.000 & 0.438 & 0.357 & 0.331 & 0.000 & 0.096 & 0.019 \\
                0.254 & 0.268 & 0.453 & 0.000 & 0.000 & 0.000 & 0.337 & 0.196 & -0.196 \\
                0.148 & 0.156 & 0.264 & 0.139 & 0.113 & 0.096 & 0.196 & 0.145 & -0.108 \\
                -0.148 & -0.156 & -0.264 & 0.028 & 0.023 & 0.019 & -0.196 & -0.108 & 0.115 \\
        \end{array}\right\rbrack
    \end{equation*}
    The test MNU is designed such that input relevance increases with decreasing input axis $m\in\mi{M}$ for $\tte[\mi{8}]{f}$. The final, exceptional output $\tte[9]{f}$ instead exhibits decreasing input relevance with decreasing $m\in\mi{M}$. As a result, only the bottom right corner of the closed Sobol' matrix shows substantial changes (10\% or more) as $\mi{m}$ decreases. The rest of the closed Sobol' matrix changes only gradually as (larger $m$) inputs of little relevance are marginalized. This is particularly true of the chosen Ishigami and Sobol' G functions: the upper-left sub-matrices $\tte[\mi{6\x 6}]{S_{\mi{5}}}$, $\tte[\mi{6\x 6}]{S_{\mi{4}}}$ and $\tte[\mi{6\x 6}]{S_{\mi{3}}}$ are very similar.
    Furthermore, additional input dimensions do not change the closed Sobol' matrix as they do not affect the test MNU:
    \begin{equation*}
        {S_{\mi{m}}}={S_{\mi{5}}} \quad \forall m\geq 5
    \end{equation*}

\begin{figure}[htbp]
        \centering
        \includegraphics[scale=0.7]{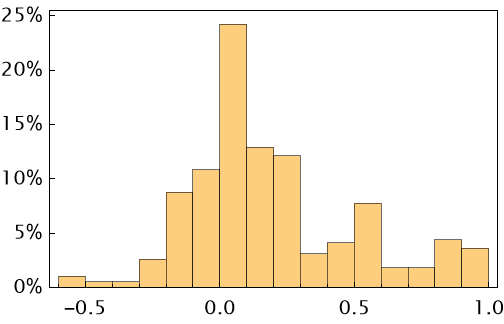}\hspace{1cm}
        \includegraphics[scale=0.7]{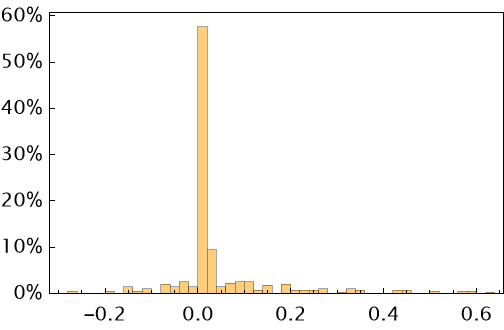}
        \caption{Frequency of closed and first order Sobol' matrix element values across 9 outputs and 5 inputs.} \label{fig:Func:S}
    \end{figure}

    Histograms of Sobol' matrix element values across the test MNU are shown in \Cref{fig:Func:S}. The benchmarking of first order Sobol' matrices reveals superb agreement, but will not be reported. We shall ignore these results because the first order Sobol' matrices of the test MNU have so few elements of appreciable magnitude that any calculation biased toward zero could produce superb agreement. Instead, benchmarking will use the closed Sobol' matrices, as these have a much more even spread of values for the test MNU.

    The most pertinent GSA features of the test MNU have just been described, and most other interesting GSA information can be inferred from them. Benchmarking consists of discovering this information as if we did not already know it. The discovered information is then compared with the ground truth just described. Effective benchmarking is stringent about this: the test MNU should be obscured in realistic ways until its Sobol' matrices are no longer discernible. Obscurity in real-world GSA almost always concerns data quantity and/or quality. 
    
    Data quantity is benchmarked by $N$, the number of samples in the design matrix being benchmarked. The $N$ samples are always arranged in a space-filling latin hypercube design \cite{Viana2015,Saltelli2008} over $M$ dimensional input space. 
    
    Data quality is impaired by adding iid (independent, identically distributed) Gaussian noise $E\rv{e}$ of magnitude $E$ to the normalized MNU, using the identity matrix $\diag[]{1}$ to define
    \begin{subequations} \label[definition]{def:Func:y}
        \begin{align}
            \te[l\x \mi{N}]{y(\tte[\mi{M\x N}]{u})} &\deq (1+E^{2})^{-1/2} \left(\te[l\x \mi{N}]{\bar{f}(\tte[\mi{M\x N}]{u})} + E \te[\mi{(l+1) N - l N}]{\rv{e}}\right) \\
            \te[l\x \mi{N}]{\bar{f}(\tte[\mi{M\x N}]{u})} &\deq \left(\te[l\x \mi{N}]{f(\tte[\mi{M\x N}]{u})}-\evt{\te[l]{f(\tte[\mi{M\x N}]{u})}}{\mi{N}}\right) / \;\devt{\te[l]{f(\tte[\mi{M\x N}]{u})}}{\mi{N}}\label[definition]{def:Func:fbar}\\
            \te[\mi{L N}]{\rv{e}} &\phantom{:}\sim \gauss{\tte[\mi{LN}]{0}}{\diag[\mi{LN\x LN}]{1}}
        \end{align}            
    \end{subequations}
    Data quality is benchmarked by $1/E$, which is the signal-to-noise ratio of $\tte[l\x \mi{N}]{y(\tte[\mi{M\x N}]{u})}$.

    Test MQUs are constructed by fitting a GP with an RBF kernel \cite{Rasmussen.Williams2005} to $\tte[l\x \mi{N}]{y(\tte[\mi{M\x N}]{u})}$. 
    The suite of test MQUs exhaustively encompasses
    \begin{align*}
        M \in \{&5,7\} \QT{;} L=9 \\
        N \in \{&30,   50,   70,   90,  110,  130,  150,  170,  190,  210,  230,
        260,  290,  320,  360,  400,  440,  480,  525,  575, \\ 
         &630,  690,
        755,  825,  900,  980, 1065, 1105, 1200, 1300, 1410, 1530, 1660,
       1800, 1950, 2110, \\
       &2280, 2460, 2710, 2930, 3170, 3430, 3710, 4000,
       4300, 4600, 4920\}\\
       E \in \{&0.0025, 0.005, 0.01, 0.025, 0.05, 0.075, 0.1, 0.2, 0.3, 0.4, 0.5, 0.6, 0.7, 0.8, 0.9, \\
       &1.0, 1.2, 1.5, 2.0, 5.0\}
    \end{align*}
    When data is very scarce ($N\leq 100$) or noisy ($E\geq 1.0$) it is ridiculously hopeful to expect good results with 5 or 7 input dimensions.

    The quality of regression is monitored by 2-fold sampling, wherein a latin hypercube of $2N$ MQU samples is randomly split in half. One half is used for regression (training) and the other for validation (testing) in the first fold, the roles reversed in the second fold. 
    
    This makes a total of nearly two thousand $L=9$ test MQUs, embodying nearly seventeen thousand individual outputs and just over 1.8 million closed Sobol' matrix elements over $M\in\set{5,7}$. Although the results are somewhat similar, $M=5$ gives a nice spread of Sobol' matrix element magnitudes, whereas $M=7$ invites the curse of dimensionality \cite{Altman.Krzywinski.2018,Binois2021} by incorporating irrelevant input dimensions. The Sobol' matrices for these irrelevant dimensions is of particular interest, as it enables dimension reduction via GSA. 
    
    The RomCom package \cite{RomCom2025} was used to implement this suite of test MQUs in a Python script to benchmark Sobol' matrix calculations and their errors.

\section{Benchmarking}\label{sec:Benchmarks}
    This Section describes the results of applying the approach developed here to the test MQUs described in the previous Section. Results are presented as heat maps of accuracy (bluer, darker hues) increasing with data quantity ($\log_{10}N$ on the x-axis) and data quality ($-\log_{10}E$ on the y-axis). Benchmarking is essentially observing the loss of accuracy and/or reliability observed as the heat map is traversed from upper right to lower left.
    A more detailed presentation, including GP calculation details, is given in \cite{MiltonRobert.etal.2025}.

\begin{figure}[b]
        \centering
        \includegraphics[scale=0.65]{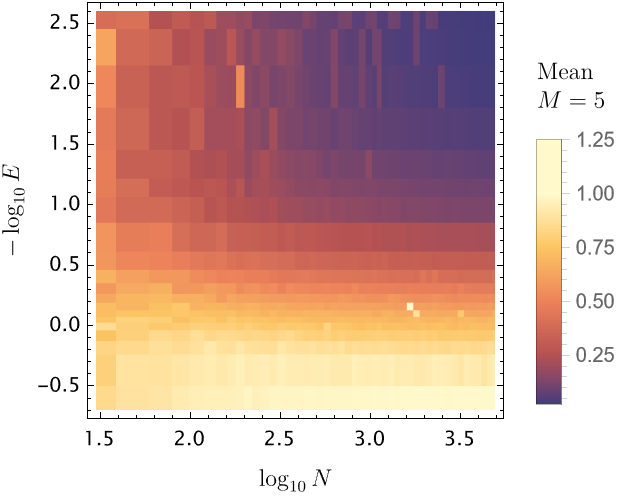}\hspace{0.5cm}
        \includegraphics[scale=0.65]{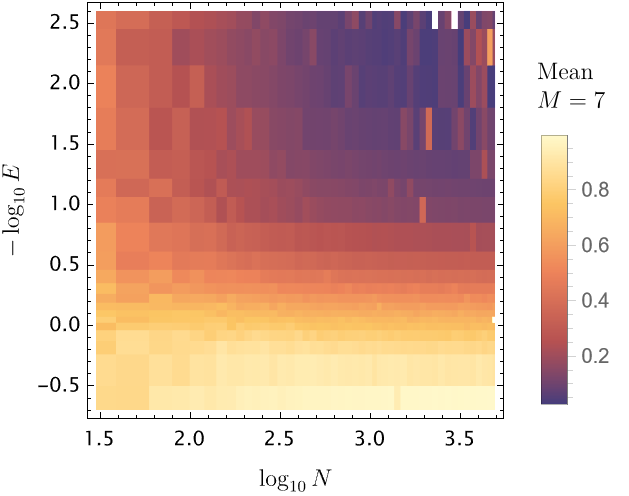}
        \caption{The Standard Deviation (SD) of prediction for the benchmark GPs. Each plot point is the mean over 2 folds and 9 outputs, so extreme outliers are visible as blemishes.} \label{fig:Bench:SD}
    \end{figure}

\begin{figure}[t]
        \centering
        \includegraphics[scale=0.65]{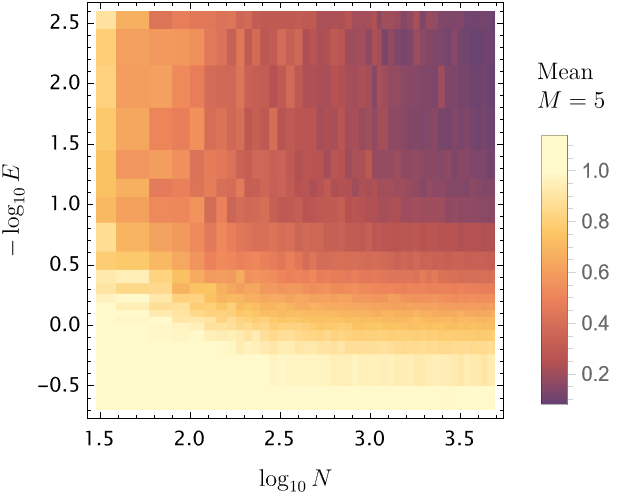}\hspace{0.5cm}
        \includegraphics[scale=0.65]{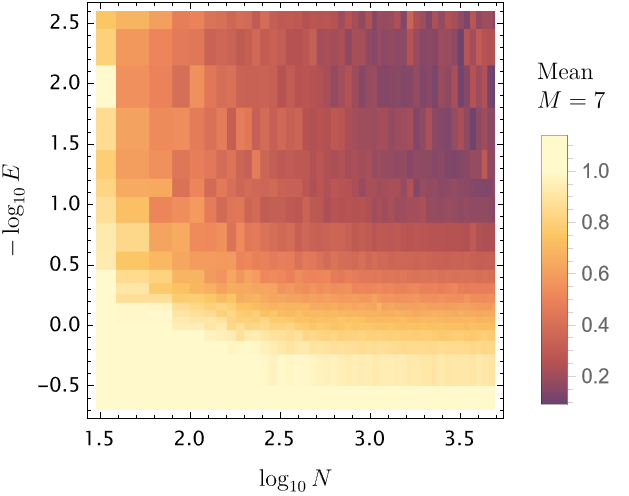}
        \caption{The Root Mean Square Error (RMSE) of prediction for the benchmark GPs, obtained from 2-fold validation. Each plot point is the mean over 2 folds and 9 outputs, so extreme outliers are visible as blemishes.} \label{fig:Bench:RMSE}
    \end{figure}

    Prior to examining the results of GSA, we consider the accuracy and reliability of the underlying MQUs, especially where data is scarce or noisy. \Cref{fig:Bench:SD} shows the standard deviation of GP predictions, \Cref{fig:Bench:RMSE} the Root Mean Square Error of GP predictions from 2-fold validation. The former indicates how accurately the GP thinks it predicts the MNU, the latter how good the predictions actually are against independent, noisy test data. The two measures show identical trends, although the SD heat maps are uniformly darker as they underestimate RMSE by about 20\%: to this degree, GP prediction is not as good as it thinks it is.
    
    \begin{figure}[b]
        \centering
        \includegraphics[scale=0.65]{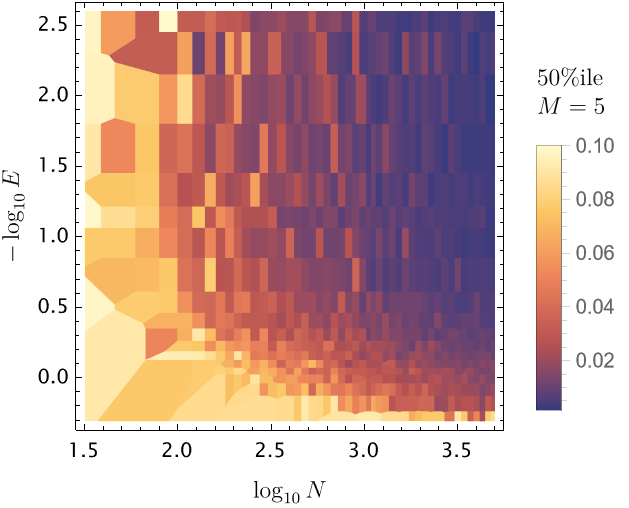}\hspace{0.5cm}
        \includegraphics[scale=0.65]{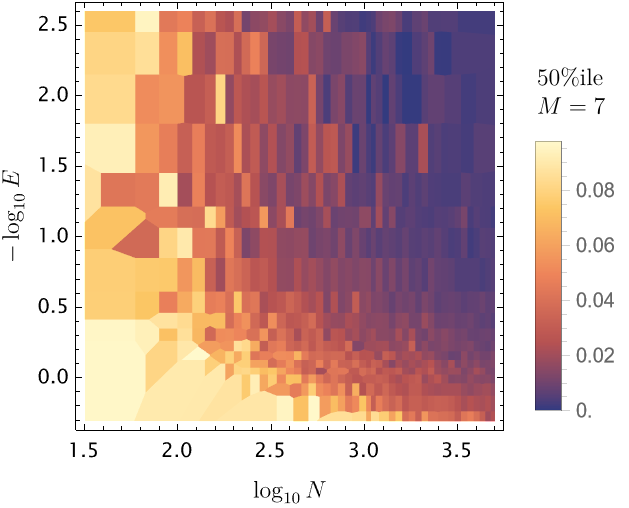}
        \vspace{0.2cm}
        \includegraphics[scale=0.65]{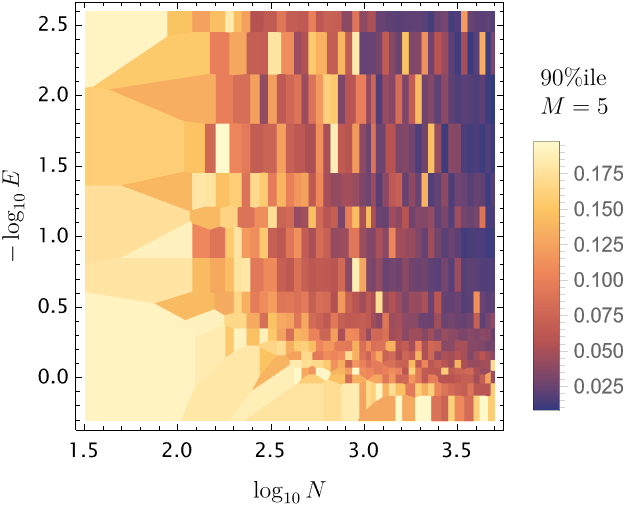}\hspace{0.5cm}
        \includegraphics[scale=0.65]{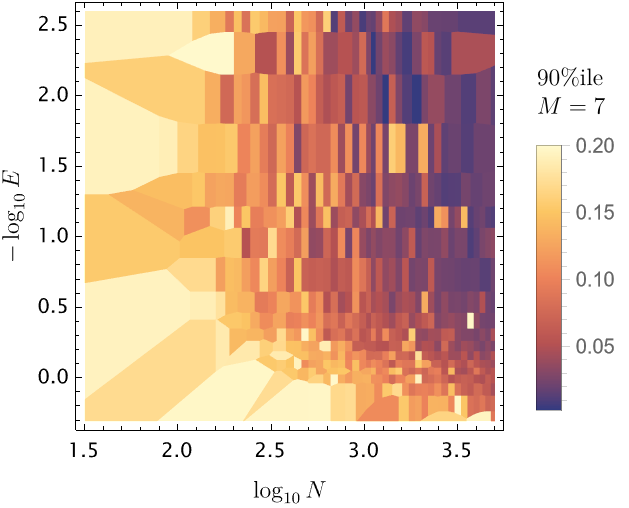}\\
        \caption{The absolute error $\tte[l\x l^{\prime}]{A_{\mi{m}}}$ of the closed Sobol' matrix element $\tte[l\x l^{\prime}]{S_{\mi{m}}}$ calculated from MQU moments. Ground truth is the MNU closed Sobol' matrix reported in \Cref{sec:Func}. Each plot point is the median (top row) or 90\% quantile (bottom row) $A$ over 2 folds, $L^{2}=81$ matrix elements and $M\in\set{5,7}$ closed matrices.The top row indicates accuracy, the bottom row reliability.}\label{fig:Bench:A}
    \end{figure}

    Recall from \Cref{def:Func:y} that each underlying MQU (GP) has been scaled to a total variation over training data of $\devt{\te[l]{y(\tte[\mi{M\x N}]{u})}}{\mi{N}}=1$. Therefore, predictions with an $\T{RMSE}=1$ are uninformative: the mean of the training data is just as good a predictor as the GP. Predictions with an $\T{RMSE}>1$ are actually worse than the overall sample mean, because the GP is overfit. The lower portion of each figure shows this occurring for $E\geq 10^{0.5}$: when noise overwhelms signal in the training data the GP is effectively useless for prediction. 
    
    There is very marked and consistent improvement (darkening) in predictive GP performance with increasing data quantity ($N$) and quality ($1/E$), and on decreasing $M$ from 7 to 5. The mean prediction error is 30\% or more whenever the noise-signal ratio $E>20\%$ or the sample size $N<100$. This description is largely valid for each GP, even though the results are displayed as mean values so that extreme outliers -- due to numerically unstable GPs -- are visible.
    Numerical instability results from extreme overfitting, and is a well understood, readily diagnosed and sometimes remediable pitfall of GP regression. It is highly susceptible to the curse of dimensionality \cite{Altman.Krzywinski.2018,Binois2021}. Given $k$-fold validation, one would normally remove these GPs, or try alternative fit optimizations to obtain a reliable GP regression.

    \begin{figure}[b]
        \centering
        \includegraphics[scale=0.65]{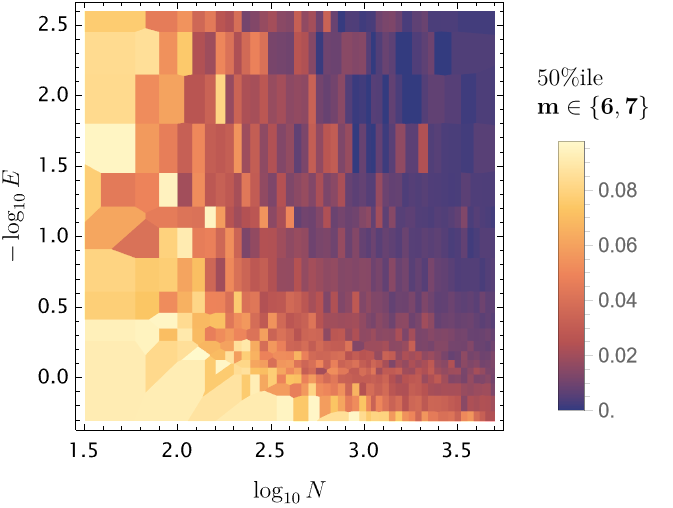}
        \includegraphics[scale=0.65]{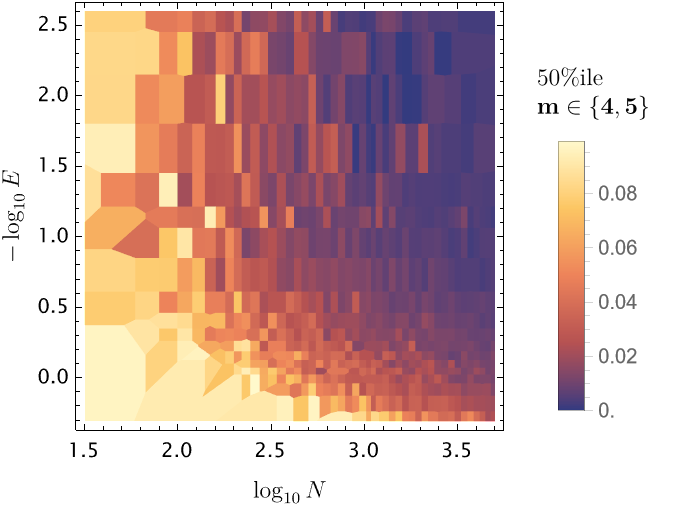}\\
        \vspace{0.2cm}
        \includegraphics[scale=0.65]{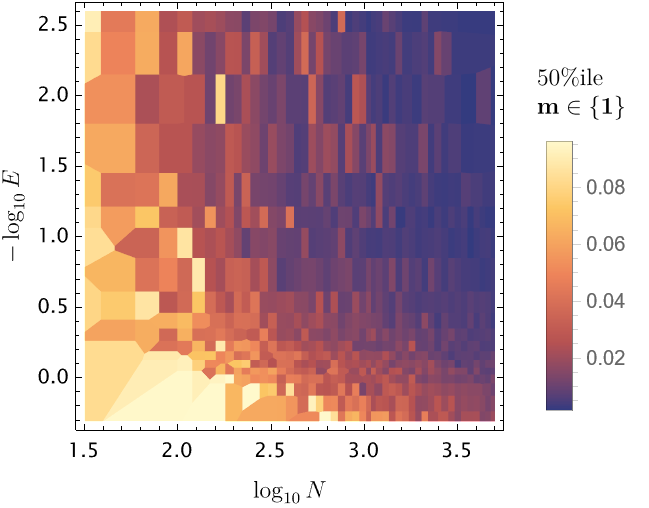}
        \includegraphics[scale=0.65]{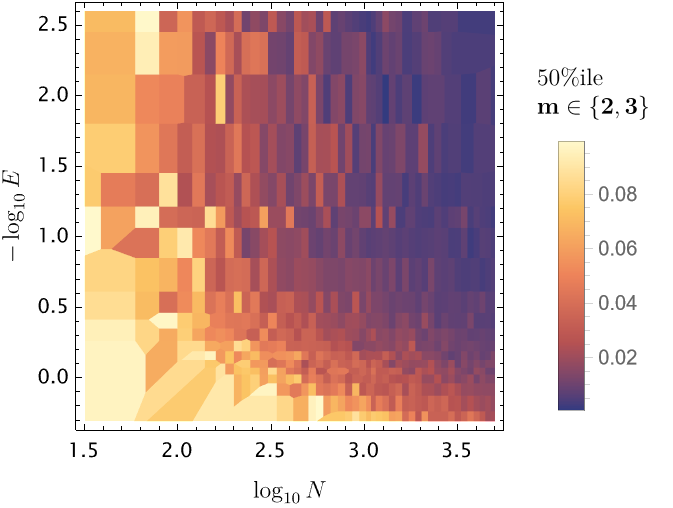}
        \caption{The absolute error $\tte[l\x l^{\prime}]{A^{T}_{\mi{M-m}}}$ of the total Sobol' matrix element $\tte[l\x l^{\prime}]{S^{T}_{\mi{M-m}}}$ calculated from MQU moments. Ground truth is $S_{\mi{M}}-S_{\mi{m}}$ calculated from the MNU closed Sobol' matrices reported in \Cref{sec:Func}. Each plot point is the median $A^{T}$ over 2 folds, $L^{2}=81$ matrix elements and one or two total matrices. Only MQUs with $M=7$ have been used.}\label{fig:Bench:Am}
    \end{figure}

    The remainder of this Section presents GSA benchmark results for Sobol' matrices calculated from the moments of the GPs just described.
    In order to assess the performance of GSA, we have not removed numerically unstable GPs using $k$-fold validation as should be done in practice. Instead, the results have been cleaned by filtering out numerically impossible GSA results. Any Sobol' index $\tte[l\x l]{S_{\mi{m}}}\not\in \unitia$ or Sobol' matrix standard error $\tte[l\x l^{\prime}]{T_{\mi{m}}}\not\in \unitia$ prompts us to remove $\tte[l]{y},\tte[l^{\prime}]{y}$ entirely from that fold $(k,M,N,E)$. If needed, such filtering should be applied in practice. On the whole, it seems that the Sobol' matrix error is a most sensitive indicator of a numerically unstable GP.

    \Cref{fig:Bench:A} shows the absolute error in Sobol' matrix elements calculated from MQU (i.e. GP) moments. The accuracy and reliability is far greater than the GPs' predictive accuracy: predictions with accuracy much worse than 30\% generate median Sobol' matrices to an accuracy of 0.05-0.10. Provided $N\geq 100$ and $E\leq 50\%$, one can confidently expect 90\% of calculated Sobol' matrices to be accurate within 0.15, and most Sobol' matrices (the median) to be good to within 0.02-0.04. $N=500$ samples nearly always provide Sobol' matrices to within 0.01 of the true value. In line with GP prediction accuracy, the estimates deteriorate with decreasing $N$ and increasing $M$, as expected. However, deterioration with increasing $E$ is much less apparent: increasing noise has little effect on GSA via MQU moments. This desirable behaviour was discussed immediately following \Cref{eq:Mom:W}, and should apply widely.

    Next we examine how the accuracy and reliability of GSA is affected by the input axes marginalized. For the purposes of dimension reduction, inputs axes $\mi{m}$ whose complementary total Sobol' matrix $S^{T}_\mi{M-m}$ is small are irrelevant, and may be ignored. The test MNU has been arranged such that input relevance mostly decreases as $m$ increases, and must be zero for $m>5$. \Cref{fig:Bench:Am} shows the median absolute error in $S^{T}_\mi{M-m}$ as input axes are removed one or two at a time. There is remarkable stability and consistency: the accuracy of the total index shows little or no dependence on $\mi{m}$, even though its value $S^{T}_\mi{M-m}$ actually follows a strong trend.

    \begin{figure}[t]
        \centering
        \includegraphics[scale=0.65]{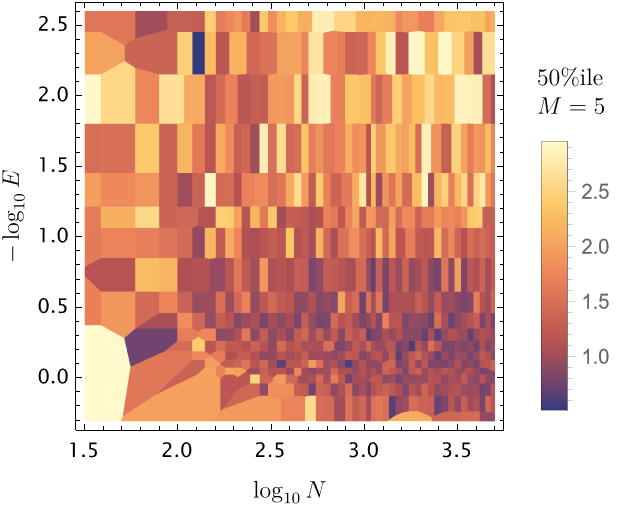}\hspace{0.5cm}
        \includegraphics[scale=0.65]{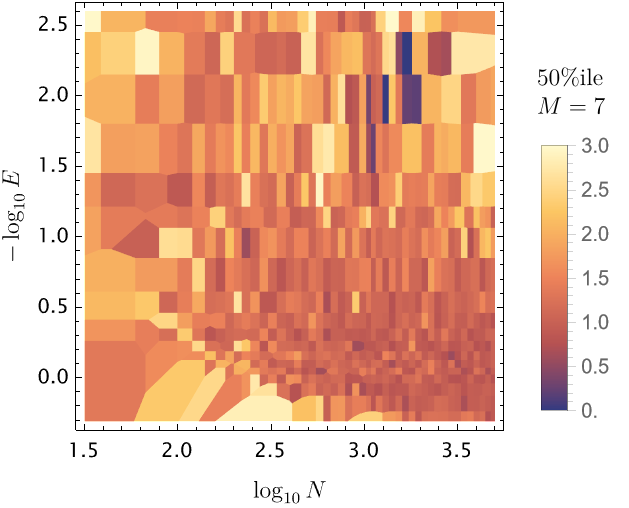}
        \vspace{0.2cm}
        \includegraphics[scale=0.65]{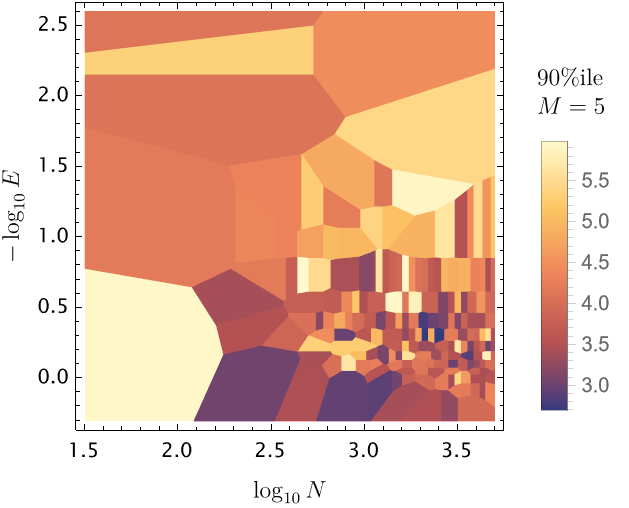}\hspace{0.5cm}
        \includegraphics[scale=0.65]{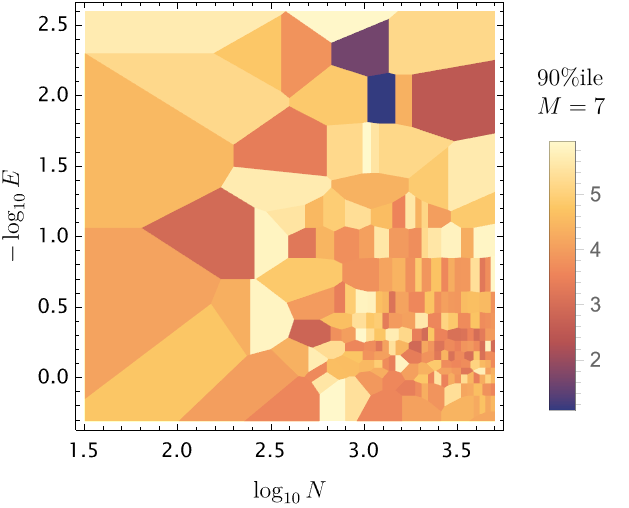}\\
        \caption{The standardized score $\tte[l\x l^{\prime}]{A_{\mi{m}}}/\tte[l\x l^{\prime}]{T_{\mi{m}}}$ of the closed Sobol' matrix element $\tte[l\x l^{\prime}]{S_{\mi{m}}}$ with standard error $\tte[l\x l^{\prime}]{T_{\mi{m}}}$ calculated from MQU moments. Ground truth is the MNU Sobol' matrix reported in \Cref{sec:Func}. Each plot point is the median (top row) or 90\% quantile (bottom row) $A/T$ over 2 folds, $L^{2}=81$ matrix elements and $M\in\set{5,7}$ closed matrices. The top row indicates accuracy, the bottom row reliabilty.}\label{fig:Bench:Z}
    \end{figure}

    MQU moments provide standard errors $\tte[l\x l^{\prime}]{T_{\mi{m}}}$ on Sobol' matrix elements in \Cref{sec:GSAS}. Normalizing the results depicted in \Cref{fig:Bench:A} by these standard errors produces the standardized score for the closed Sobol' matrix element $\tte[l\x l^{\prime}]{S_{\mi{m}}}$ shown in \Cref{fig:Bench:Z}. One does not expect $\tte[l\x l^{\prime}]{\rv{S}_{\mi{m}}}$ to be normally distributed, so standardized scores are not expected to follow rules such as 95\% of results score 2 or less. Rather, one is looking for stable and consistent scores over any choice of $E,N,M$. This means that the standard error $\tte[l\x l^{\prime}]{T_{\mi{m}}}$ has a meaning which is largely independent of context, and provides a consistent guide to the accuracy and reliability of $\tte[l\x l^{\prime}]{S_{\mi{m}}}$. Looking at \Cref{fig:Bench:Z}, one sees a median standardized score consistently less than 3. The 90\% quantile is uniformly less than 6. However, the 90\% quantile shows a marked increase with increased input dimensionality from $M=5$ to $M=7$. The reliability -- but not necessarily the accuracy -- of standard errors on Sobol' matrices computed from GP moments are very susceptible to the curse of dimensionality.

\section{Conclusions}\label{sec:Conc}
    In this paper we have presented the theory and implementation of Sobol' indices extended to multi-outputs in the form of Sobol' matrices. These essentially ascribe the correlation between output components to the influence of various inputs. Every element of a Sobol matrix lies within $\lbrack-1,1\rbrack$. The diagonal elements of the Sobol' matrix are precisely the Sobol' indices of the output components, and the closed Sobol' matrix $\rv{S}_{\mi{M}}$ considering all inputs is just the multi-output correlation matrix. This is deliberately intended to provide insights into the influence of inputs on the linkages (correlations) between output components.
    
    The formal development assumes a multi-output model with quantified uncertainty (MQU) which is essentially a stochastic process (SP). The Sobol' matrix and its standard error have been expressed in terms of the first two moments of this SP. When benchmarked using GPs regressed to standard test functions with noise, the calculation is remarkably accurate up to $M=7$ input dimensions. Sobol' matrices of 0.02 accuracy are usually obtained with just 300 samples, and even 90 samples will usually provide estimates accurate to within 0.08 of the true value. It is surely advantageous that Sobol' matrices and their uncertainties are not directly affected by homoskedastic uncertainty or noise in an MQU. However, when the noise to signal ratio exceeds $10^{0.5}$ accuracy is indirectly impaired by the woeful performance of the GP as a predictor. The accuracy of total Sobol' matrices also appears hardly affected by the number of relevant or irrelevant inputs included up $M=7$. The standard error of Sobol' matrices indicate that most Sobol' matrices are accurate to within 3 standard errors, and 90\% are accurate to within 6 standard errors. It should be ascertained in future whether these observations persist with more input dimensions $M>7$. This would be very useful for reducing high dimensional systems to their relevant inputs.

    Clearly the GP surrogate MQU moment method of GSA developed here looks very powerful when input dimensionality is moderate and data is scarce, or noisy, or both. It does not scale well with the number of outputs or samples: if data is abundant, other methods may be preferable. In any case, if one wishes to study the influence of various inputs on the correlations between outputs, Sobol' matrices provide a coherent framework for doing so.


\bibliographystyle{amsplain}
\bibliography{SobolMatrices}

\end{document}